\documentclass{amsart}
\usepackage{amsmath, amssymb}
\newtheorem{thm}{Theorem}[section] 
\newtheorem{cor}[thm]{Corollary}
\newtheorem{lem}[thm]{Lemma}
\newtheorem{defn}[thm]{Definition}
\newtheorem{prop}[thm]{Proposition}

\newtheorem{preremark}[thm]{Remark}
\newenvironment{remark}%
  {\begin{preremark}\upshape}{\end{preremark}}
\newtheorem{preexample}[thm]{Example}
  {\begin{preexample}\upshape}{\end{preexample}}

\numberwithin{equation}{section}

\numberwithin{equation}{section}

\newcommand{\ten}{\otimes}

\newcommand{\w}{\mathbf{w}}
\newcommand{\abs}[1]{\lvert#1\rvert}
\newcommand{\comS}[1]{[#1]_S}
\newcommand{\inv}{^{-1}}
\newcommand{\MapzV}{\Map_{z_1,z_2,\dots,z_n}(V^{\ten n})}

\newcommand{\nopS}[1]{:\!#1\!:_{S}}
\newcommand{\Ratz}{\Rat_{z_1,z_2,\dots,z_n}}
\newcommand{\subz}{_{z_1,\dots,z_n}}

\DeclareMathOperator{\End}{End}
\DeclareMathOperator{\GL}{GL}
\DeclareMathOperator{\Hol}{Hol}
\DeclareMathOperator{\Hom}{Hom}
\DeclareMathOperator{\Map}{Map}
\DeclareMathOperator{\Rat}{Rat}
\DeclareMathOperator{\Res}{Res}
\DeclareMathOperator{\Sing}{Sing}

\begin{document}

\title[$H_D$-Quantum Vertex Algebras]{$H_D$-Quantum Vertex
  Algebras and Bicharacters}
\author{Iana I. Anguelova}
\author{Maarten J. Bergvelt}

\address{Anguelova:
Centre de Recherches Mathematiques (CRM)\\
Universit\'e de Montreal\\
Montreal,  Quebec H3C 3J7
}
\email{anguelov@crm.umontreal.ca}
\address{Bergvelt: Department of Mathematics\\ University of Illinois\\
  Urbana-Champaign\\ Illinois 61801}
\email{bergv@uiuc.edu}


\begin{abstract}
  We define a new class of quantum vertex algebras, based on the Hopf
  algebra $H_D=\mathbb{C}[D]$ of "infinitesimal translations"
  generated by $D$. Besides the braiding map describing the
  obstruction to commutativity of products of vertex operators,
  $H_D$-quantum vertex algebras have as main new ingredient a
  "translation map" that describes the obstruction of vertex operators
  to satisfying translation covariance. The translation map also
  appears as obstruction to the state-field correspondence being a
  homomorphism.

  We use a bicharacter construction of Borcherds to construct a large
  class of $H_D$-quantum vertex algebras. One particular example of
  this construction yields a quantum vertex algebra that contains the
  quantum vertex operators introduced by Jing in the theory of
  Hall-Littlewood polynomials.
\end{abstract}

\maketitle
\tableofcontents

\section{Introduction}
\label{sec:intro}

Vertex operators were introduced in the earliest days of string theory
and play now an important role of such areas of mathematics as
representation theory, algebraic topology and random matrices. Vertex
algebras were introduced to axiomatize the properties of vertex
operators.

Similarly, quantum vertex operators were discovered in integrable
models in statistical mechanics and in connection with theory of
symmetric polynomials and the theory of quantum affine algebras. One
would like to have theory of quantum vertex algebras to axiomatize the
properties of quantum vertex operators. In this paper we introduce and
study a class of quantum vertex algebras that produce the quantum
vertex operators related to Hall-Littlewood polynomials.

Recall that a vertex operator on a space $V$ is a series
$a(z)=\sum_{n\in \mathbb{Z}}a_{(n)}z^{-n-1}$, $a_{(n)}\in \End(V)$,
satisfying some extra conditions. We call vertex operators
\emph{local} with respect to each other if the commutator is a sum of
derivatives of delta distributions:
\[
[a(z_1),b(z_2)]=\sum_{n=0}^N c_n(z_2)\partial_{z_2}\delta(z_1,z_2).
\]
For \emph{quantum vertex operators} this will not longer be true: one
needs a \emph{braiding} map $S^{(\tau)}_{z_1,z_2}\colon
b(z_2)a(z_1)\mapsto ba(z_2,z_1)$, where $ba(z_2,z_1)$ is some other
$\End(V)$-valued series. Then we should have $S$-locality,
\cite{MR2002i:17022}, i.e., we need that the \emph{braided commutator}
\[
\comS{a(z),b(z_2)}=a(z_1)b(z_2)-ba(z_2,z_1)
\]
is a sum of derivatives of delta distributions. One of the goals of
this papers is give explicit examples of quantum vertex algebras where
one can easily calculate both the quantum vertex operators and their
braiding.

There are several proposals for what a quantum vertex algebra should
be. There is Borcherds' theory of (A, H, S)-vertex algebras, see
\cite{MR1865087}, the Etingof-Kazhdan theory of quantum vertex
operator algebras, \cite{MR2002i:17022} , and the Frenkel-Reshet'ikin
theory of deformed chiral algebras, see \cite{q-alg/9706023}. (H. Li
has developed the Etingof-Kazhdan theory further, see for example
\cite{MR2220654}, \cite{MR2215259}.)

The Borcherds theory is based on the observation that products and
iterates of vertex operators in vertex algebras are expansions of
rational functions in multiple variables. The idea then is to start
with these rational objects instead of constructing them after the
fact from the vertex operators. Instead of a single vector space $V$
on which the vertex operators act, one has for any integer $n\ge 1$
the space $V(n)$ of ''rational vertex operators'' in $n$ variables.
This is quite beautiful idea, and is easily adapted to include quantum
vertex algebras of (A, H, S)-type. However, even for classical vertex
algebras it seems not known how to include such basic examples as
affine vertex algebras in the (A, H, S)-framework. In this paper we
therefore we prefer to develop a theory that is closer to the usual
theory, with a single underlying vector space $V$. We do take,
however, from Borcherds' paper the idea of a \emph{bicharacter} as a
method to construct examples: we will use bicharacters both to produce
the vertex operators and the braiding. (See \cite{Anguelova:thesis}
for more details.)

The Etingof-Kazhdan theory is very close to the classical theory, in
fact so close that it is not suitable to describe quantum vertex
operators related to symmetric polynomials. Briefly, in the usual
theory (and in \cite{MR2002i:17022}) vertex operators $a(z)$ satisfy
translation covariance of the form
\begin{equation}
e^{\gamma D}a(z_1)e^{-\gamma D}=a(z_1+\gamma),\label{eq:transcovariance}
\end{equation}
where $D\colon V\to V$ is the infinitesimal translation operator and
we expand in positive powers of $\gamma$. If we introduce notation
$a(z)b=Y_z(a\otimes b)$ we can write this, since
$\partial_za(z)=(Da)(z)$, as
\[
e^{\gamma D} Y_z=Y_z\circ(e^{\gamma D}\otimes e^{\gamma D}),
\]
making clear the similarity of a vertex algebra with a associative
ring $M$ with a group action (where $g m(a\ten b)=m(ga\ten g b)$ if
$m$ is the multiplication in $M$, and $g\in G$, $a,b\in M$, see
Appendix \ref{sec:BraidAlgSym}).

It was shown in the thesis \cite{Anguelova:thesis} that
\eqref{eq:transcovariance} can not hold in the case of quantum vertex
operators related to symmetric polynomials. Also, not unrelated, the
braiding map $S_{z_1,z_2}$ in \cite{MR2002i:17022} is in fact assumed
to be of the form $S_{z_1,z_2}=\tilde S_{z_1-z_2}$, where $\tilde S_z$
is a function of a single variable. It is also shown in
\cite{Anguelova:thesis} that this does not holds for symmetric
polynomials.

In this paper we introduce the notion of an $H_D$-quantum vertex
algebra (where $H_D=\mathbb{C}[D]$ is the Hopf algebra of
infinitesimal translations), generalizing \cite{MR2002i:17022} in
various ways. First we need to relax the translations covariance
\eqref{eq:transcovariance}. We introduce, besides the braiding map
$S^{(\tau)}_{z_1,z_2}$, also another map $S^{(\gamma)}_{z_1,z_2}$ on
$V\otimes V$ such that we get instead
\[
e^{\gamma D}Y_{z_1}\circ S^{(\gamma)}_{z_1,z_2}=Y_{z_1}\circ
(e^{\gamma D}\otimes e^{\gamma D}).
\]
Both $S^{(\gamma)}_{z_1,z_2}$ and $S^{(\tau)}_{z_1,z_2}$ are rational
functions of both $z_1$ and $z_2$, not just of the difference
$z_1-z_2$ as in \cite{MR2002i:17022}. Another difference is that in
\cite{MR2002i:17022} vertex operators satisfy a braided version of
skew-symmetry:
\begin{equation}
  \label{eq:skewsymY}
Y_z\circ S^{(\tau)}_{z,0}(a\otimes b)= e^{zD}Y(b,-z)a.
\end{equation}
This relation does not make sense for quantum vertex operators coming
from symmetric polynomials: the braiding $S^{(\tau)}_{z_1,z_2}$ is in
general \emph{singular} for $z_2=0$. This motivates us to take as
basic building block of the theory not the vertex operator $Y_z$, but
the \emph{two-variable} vertex operators $X_{z_1,z_2}\colon V\otimes
V\to V[[z_1,z_2]][z_1\inv,(z_1-z_2)\inv][[t]]$. 
We can define then $Y$ by $Y(a,z)b=X_{z,0}(a\ten b)$, but $Y$ does not
longer satisfy \eqref{eq:skewsymY}. See Corollary \ref{cor:skewsym}
for the version of skew-symmetry that holds for $H_D$-quantum vertex
algebras. 

Conversely, if we start with $Y$,
we can introduce $X_{z_1,z_2}$ by analytic continuation: we have the
expansion
\begin{equation}
  \label{eq:expanXgivesY}
i_{z_1;z_2}X_{z_1,z_2}(a\otimes b)=Y(a,z_1)Y(b,z_2)1,
\end{equation}
where $i_{z;w}$ is the expansion in the region $\abs z>\abs
w$. See Section \ref{sec:AltAx} for details and an alternative
definition of $X_{z_1,z_2}$.

Note that for a classical vertex algebra (and also for an
Etingof-Kazhdan quantum vertex operator algebra) the translation map
$S^{(\gamma)}_{z_1,z_2}$ is the identity, so that in this case
\[
X_{z_1,z_2}(a\ten b)=e^{z_2}Y(a,z_1-z_2)b\in
V[[z_1,z_2]][(z_1-z_2)\inv].
\]
In particular in these cases $X_{z_1,z_2}(a\ten b)$ is not singular
for $z_1=0$. We consider a more general theory where in
$X_{z_1,z_2}(a\ten b)$ poles in $z_1$ are allowed (and in fact are
necessary to be able to treat the quantum vertex operators associated
with the Hall-Littlewood polynomials).

In the construction of quantum vertex algebras one or more quantum
parameters will appear. They can usually be thought of as describing
the deformation away from an ordinary vertex algebra. We should
mention that, just as when quantizing universal enveloping algebras,
there are two ways of interpreting the quantum parameters in quantum
vertex algebras. Either the quantum parameters are independent formal
variables or they are complex numbers.  The theory of Etingof-Kazhdan
follows the first approach, as opposed to the Frenkel-Reshetikhin
definition of deformed chiral algebras, which considers the
deformation parameter(s) to be complex number(s). In this paper we
also follow the first approach: we have an independent variable $t$
and a $H_D$-quantum vertex algebra $V$ is a (free) module over the
ring $\mathbb{C}[[t]]$ of formal power series in $t$. When putting
$t=0$ one gets in examples generally an ordinary vertex algebra,
although we did not require this in our axioms. Note that Li in
\cite{MR2220654}, for instance, studies a form of the Etingof-Kazhdan
axioms where the quantum parameter is a complex number.

Maybe the most important difference between our $H_D$-quantum vertex
algebras and classical vertex algebras (and the theory of
Etingof-Kazhdan) is the following. We can define a products of states:
$a_{(n)}b$, and for fields: $a(z)_{(n)}b(z)$, but is is not longer
true that that the state-field correspondence $a\mapsto Y(a,z)$ a
homomorphism of products: in general $Y(a_{(n)}b,z)\ne
Y(a,z)_{(n)}Y(b,z)$, 
see Theorem \ref{thm:StfieldCorr} for an exact
statement.

The outline of the paper is as follows.  There are three parts. In the
first part we define $H_D$-quantum vertex algebras in section
\ref{sec:H_DQuantumVertex} and study their properties in the following
sections. We derive in Section \ref{sec:braid-Jac-identity} and
\ref{sec:BrBoRSlov} fundamental identities in our quantum vertex
algebras: the braided Jacobi identity and the braided Borcherds
identity. These are used Sections \ref{sec:Scomloc},
\ref{sec:n-productsstates} and \ref{sec:NormOrdProd} to study the
$S$-commutator, $(n)$-products of states and of fields and normal
ordered products. In Section \ref{sec:WeakAssoc} we derive a weak
associativity relation. In the next part of the paper we assume that
our underlying vector space $V$ is a commutative and cocommutative
Hopf algebra, which allows us to defines bicharacters on $V$ in
section \ref{sec:BiChar}. Using bicharacters we construct a class
$H_D$-quantum vertex algebras in Section \ref{sec:QuantVertBichar} and
in the rest of this section we explore some of the properties of
bicharacter $H_D$-quantum vertex algebras. 
 In last part of the paper,
Section \ref{sec:mainexample} and the following sections, we study in
detail a single example of a Hopf algebra $V$ with a fixed bicharacter
on it. The resulting $H_D$-quantum vertex algebra is a deformation of
the familiar lattice vertex algebra based on the lattice
$L=\mathbb{Z}$ with pairing $(m,n)\mapsto mn$. Some of quantum vertex
operators in this example were used by Jing, see \cite{MR1112626}, to
study Hall-Littlewood symmetric polynomials.  In the Appendix
\ref{sec:BraidAlgSym} we describe the ``nonsingular" analog of
$H_D$-quantum vertex algebras: braided algebras with group action. In
Appendix \ref{app:braiding} we discuss the construction of braiding
maps for $H_D$-quantum vertex algebras.

\section{The Hopf Algebra $H_D$}
\label{sec:HopfH_D}

Let $H_D=\mathbb C[D]$ be the universal enveloping algebra of the
1-dimensional Lie algebra generated by $D$. $H_D$ is a Hopf-algebra,
with coproduct $\Delta_{H_D}\colon D\mapsto D\otimes 1+1\otimes D$,
antipode $S\colon D\mapsto -D$ and counit $\epsilon_{H_D}\colon
D\mapsto 0$.  $H_D$ is a fundamental ingredient in the construction of
vertex algebras, where it appears as the symmetry algebra of
infinitesimal translations in physical space. In this paper the full
Hopf algebra structure of $H_D$ will play only an explicit role when
we discuss bicharacter constructions, in the definition of
$H_D$-quantum vertex algebras in the next section only the algebra
structure will be used. However, from Borcherds' papers
\cite{MR1653021} and \cite{MR1865087} it will be clear that in fact
the Hopf algebra $H_D$ underlies the whole theory of vertex algebras
(and their quantum versions).

\section{$H_D$-Quantum Vertex Algebras}
\label{sec:H_DQuantumVertex}

Let $t$ be a variable. We will use $t$ to describe quantum
deformations, the classical limit corresponding to $t\to 0$. Let
$k=\mathbb C[[t]]$ and let $V$ be an $H_D$-module and free
$k$-module. Denote by $V[[t]]$ the space of (in general infinite) sums
\[
v(t)=\sum_{i=0}^\infty v_i t^i,\quad v_i\in V.
\]
In case $v(t)\in V[[t]]$ has only finitely many nonzero terms we can
identity it with an element of $V$. In the same way will consider
spaces such as $V[[z]][z\inv][[t]]$ consisting of sums
\[
v(z,t)=\sum_{i=0}^\infty v_i(z) t^i,\quad v_i\in V[[z]][z\inv].
\]
We will also consider rational expressions in multiple variables and
their expansions. For instance for a rational function in $z_1$, $z_2$
with only possibly poles at $z_1=0$, $z_2=0$ or $z_1-z_2=0$ we can
define expansion maps
\begin{align*}
  i_{z_1;z_2} \colon \frac1{z_1-z_2}&\mapsto
  \sum_{n\ge0}z_1^{-n-1}z_2^n, &\quad \frac 1{z_1}&\mapsto \frac
  1{z_1}, &\quad \frac 1{z_2}&\mapsto \frac 1{z_2},\\
  i_{z_2;z_1} \colon \frac1{z_1-z_2}&\mapsto
  -\sum_{n\ge0}z_2^{-n-1}z_1^n,&\quad \frac 1{z_1}&\mapsto \frac
  1{z_1},&\quad  \frac 1{z_2}&\mapsto \frac 1{z_2},\\
  i_{z_2;z_1-z_2} \colon \frac1{z_1} &\mapsto
  \sum_{n\ge0}z_2^{-n-1}(z_1-z_2)^n,&\quad \frac 1{z_2}&\mapsto \frac
  1{z_2}, &\quad \frac 1{z_1-z_2}&\mapsto \frac 1{z_1-z_2}.
\end{align*}
We will write $i_{z_1,z_2;w_1}$ for $i_{z_1;w_1}i_{z_2;w_1}$, and
$i_{z_1,z_2;w_1,w_2}$ for $i_{z_1,z_2;w_1}\allowbreak
i_{z_1,z_2;w_2}$. We define $i_{z_1;z_2;\dots;z_n}$ to be the
expansion in the region $\abs{z_1}>\abs{z_2}>\dots>\abs{z_n}$.

If $A\in V\otimes V$ then we define for instance
$A^{23}, A^{13}\in V^{\otimes 3}$ by $A^{23}=1\otimes A$, and
$A^{13}=a^\prime \otimes 1\otimes a^{\prime\prime}$, if
$A=a^\prime\otimes a^{\prime\prime}$.

Now we are ready to define the central concept of this paper. The
definition is rather complicated, and in Appendix
\ref{sec:BraidAlgSym} we explain a simpler version of this notion,
called a braided ring with symmetry, where the multiplication is
nonsingular.
\begin{defn}\label{defn:h_d-quantum-vertex-alg}
  Let $V$ be a free $k=\mathbb{C}[[t]]$-module and an $H_D$-module.
  An $H_D$-quantum vertex algebra structure on $V$ consists of
  \begin{itemize}
  \item $1\in V$, the vacuum vector.
  \item a (singular) multiplication map
\[
X_{z_1,z_2}\colon V^{\otimes 2}\to V[[z_1,z_2]][z_1\inv,
(z_1-z_2)\inv][[t]].
\]
\item A braiding map $S^{(\tau)}$ and a translation map $S^{(\gamma)}$
  of the form
  \begin{align*}
    S^{(\tau)}_{z_1,z_2}&\colon V^{\otimes 2}\to V^{\otimes2}[z_1^{\pm
      1},z_2^{\pm 1},(z_1-z_2)^{- 1}][[t]],\\
    S^{(\gamma)}_{z_1,z_2}&\colon V^{\otimes 2}\to
    V^{\otimes2}[z_1^{\pm1},z_2,(z_1+\gamma)^{\pm1},(z_2+\gamma),
    (z_1-z_2)^{-1} ][[t]] .
  \end{align*}
  \end{itemize}
These objects satisfy the following axioms:
\begin{description}
\item[(Vacuum)] For $i=1,2$
  \begin{align}
    X_{z_1,z_2}(a\ten 1)&=e^{z_1D}a,&     X_{z_1,z_2}(1\ten a)&=e^{z_2D}a, \label{eq:vacuumX}\\
S_{z_1,z_2}(a\ten 1)&=a\ten 1,& S_{z_1,z_2}(1\ten a)&=1\ten a. \label{eq:vacuumS}
  \end{align}
Here and below we write generically $S$ for both $S^{(\tau)}$ and
$S^{(\gamma)}$. 

\item[($H_D$-covariance)]
\begin{align}
    X_{z_1,z_2}(a\otimes Db)&=\partial_{z_2}X_{z_1,z_2}(a\otimes b), \label{eq:H_DcovX}\\
    (1\otimes e^{\gamma
      D})i_{z_1-z_2, z_2;\gamma}S_{z_1,z_2+\gamma}&=S_{z_1,z_2}(1\otimes
    e^{\gamma D}), \label{eq:H_DcovS} \\
    e^{\gamma D}X_{z_1,z_2}S^{(\gamma)}_{z_1,z_2} &= X_{z_1+\gamma,z_2+\gamma}\label{eq:H_DcovMult}.
  \end{align}

\item[(Yang-Baxter)] 
  \begin{equation}\label{eq:YBaxiom}
S_{z_1,z_2}^{12}S_{z_1,z_3}^{13}S_{z_2,z_3}^{23}=S_{z_2,z_3}^{23}S_{z_1,z_3}^{13}S_{z_1,z_2}^{12}.
  \end{equation}
\item[(Compatibility with Multiplication)] 
  \begin{align}
    S_{z_1,z_2}(X_{w_1,w_2}\otimes 1)
    &=(X_{w_1,w_2}\otimes1)i_{z_1, z_1-z_2;w_1,w_2}
    S_{z_1+w_1,z_2}^{23}S_{z_1+w_2,z_2}^{13},\label{eq:CompatXx1}\\
    S_{z_1,z_2}(1\otimes X_{w_1,w_2})
    &=(1\otimes X_{w_1,w_2})i_{z_1-z_2, z_2;w_1,w_2}
    S_{z_1,z_2+w_1}^{12}S_{z_1,z_2+w_2}^{13}.\label{eq:Compat1xX}
  \end{align}
\item[(Group Properties)] 
  \begin{align}
    S_{z_1,z_2}^{(\tau)}\circ \tau\circ S_{z_2,z_1}^{(\tau)}\circ \tau&=1_{V^{\otimes 2}},\label{eq:GrpSStau}\\
    S^{(\gamma_1)}_{z_1,z_2}S^{(\gamma_2)}_{z_1+\gamma_1,z_2+\gamma_1}&=S^{(\gamma_1+\gamma_2)}_{z_1,z_2},
    \label{eq:GrpSgam1gam2}\\
    S^{(\gamma=0)}_{z_1,z_2}&=1_{V^{\otimes 2}}.\label{eq:SatZero}
  \end{align}
\item[(Locality)] For all $a,b\in V$ and $k\ge 0$ there is $N\ge0$
  such that for all $c\in V$
\begin{multline} 
  (z_1-z_2)^NX_{z_1,0}(1\ten
  X_{z_2,0})(a\ten b\ten c)\equiv\\
  \equiv (z_1-z_2)^N X_{z_2,0}(1\ten X_{z_1,0})\left(i_{z_2
      ;z_1}S^{(\tau)}_{z_2,z_1}(b\ten a)\ten c\right)\mod
  t^k.\label{eq:localityAx}
\end{multline}

\end{description}
\end{defn}
\begin{remark}
  In the above definition $z_1,z_2,w_1,w_2,\gamma$ are independent
  commuting variables. In general one should be careful with
  specializing these variables. For instance, we can evaluate
  $X_{z_1,z_2}$ at $z_2=0$ but not at $z_1=0$, in general. For this
  reason one can not put $\gamma=-z_1$ in \eqref{eq:H_DcovMult}.\qed
\end{remark}

\begin{remark}
  The vacuum axioms \eqref{eq:vacuumX} for $z_2=0$ are
\[
X_{z_1,0}(a\ten 1)=e^{z_1D}a,\quad X_{z_1,0}(1\ten a)=a.
\]
In the literature on vertex algebras the first equation is called the
\emph{creation axiom}, and the second the vacuum axiom. In our
formalism it seems unnatural to give different names to very similar
statements, so we call in \eqref{eq:vacuumX} both vacuum axioms, as
they involve the vacuum vector 1.
\end{remark}

\section{Intermezzo on Expansions}
\label{sec:intermezzoexp}

Let $W$ be a vector space and $A(z_1,z_2)\in W((z_1))((z_2))$.  It is
well known\footnote{ See for instance the notion of compatible fields
  in Definition 7.3 of \cite{MR1653021}, \cite{math.QA/9904104}, and
  the reformulation of compatibility in \cite{MR1966260},
  \cite{MR2218823}.} that if there is an $N\ge0$ such that
\begin{equation}
  \label{eq:compatibleA}
  A_N=(z_1-z_2)^NA(z_1,z_2)\in W[[z_1,z_2]][z_1\inv,z_2\inv],
\end{equation}
then $A(z_1,z_2)$ is in the image of the (injective) map $i_{z_1;z_2}$,
i.e., there is a (unique) $X(z_1,z_2)\in
W[[z_1,z_2]][z_1\inv,z_2\inv,(z_1-z_2)\inv]$ such that we have the expansion
\begin{equation}
  A(z_1,z_2)=i_{z_1;z_2} X(z_1,z_2).
  \label{eq:AexpansX}
  \end{equation}
  In fact, we can take $X(z_1,z_2)=(z_1-z_2)^{-N}A_N$. In this case we
  have also
  \begin{equation}
    \label{eq:pnApnX}
(z_1-z_2)^N A(z_1,z_2)=(z_1-z_2)^N X(z_1,z_2).
  \end{equation}
  (Note that although $A_N$ depends on $N$, we obtain the same $X$ for
  all $N$ that make \eqref{eq:compatibleA} true , cf.
  \cite{MR1966260}.)

One way to check \eqref{eq:compatibleA} is by finding $B(z_2,z_1)\in
W((z_2))(z_1))$ such that
\begin{equation}
  \label{eq:sortoflocality}
  (z_1-z_2)^NA(z_1,z_2)=(z_1-z_2)^NB(z_2,z_1).
\end{equation}
Indeed, the LHS shows that \eqref{eq:sortoflocality} has at worst a
finite order pole in $z_2$ (by assumption on $A(z_1,z_2)$) and the RHS
that at worst it has a finite order pole in $z_1$ (by assumption on
$B(z_2,z_1)$). This means that \eqref{eq:sortoflocality} belongs to
$W[[z_1,z_2]][z_1\inv,z_2\inv] $, as we wanted to show.  In this case
we have not only that $A$ is the expansion \eqref{eq:AexpansX} of $X$,
but also that $B$ is the ``opposite'' expansion:
\[
B(z_1,z_2)=i_{z_2;z_1} X(z_1,z_2).
\]
There are generalizations to more variables $z_1, z_2,\dots,
z_n$, and to various expansion maps.

We will need slight refinements of these phenomena in case there is a
quantum parameter $t$ present. For example: 
 
\begin{lem}\label{lem:z1-z2Arational}
    Let
\[ 
A(z_1,z_2;t)\in W((z_1))((z_2))[[t]],
\]
and supppose that for all $k\ge0$ there is an $N\ge0$ such that
\begin{equation}
  \label{eq:rationality}
A^k_N\equiv (z_1-z_2)^NA(z_1,z_2;t)\mod t^k\in
W[[z_1,z_2]][z_1\inv,z_2\inv][[t]]/\langle t^k\rangle.
\end{equation}
Then there is a $X(z_1, z_2)\in
W[[z_1,z_2]][z_1\inv,z_2\inv,(z_1-z_2)\inv][[t]]$ such that
\begin{equation}
i_{z_1;z_2}X(z_1,z_2)=A(z_1,z_2).
\label{eq:expansionXisA}
\end{equation}
\end{lem}
\begin{proof}
  If \eqref{eq:rationality} holds for some $N$ we can define 
\[
X^k=(z_1-z_2)^{-N} A^k_N\in W[[z_1,z_2]][z_1\inv,z_2\inv,(z_1-z_2)\inv][[t]]/\langle t^k\rangle,
\]
and we have 
\[
i_{z_1;z_2}X^k(z_1,z_2)= A^k_N(z_1, z_2;t).
\]
Then the $X^k$s fit together to define a (unique)\\
$X(z_1, z_2;t)\in W[[z_1,z_2]][z_1\inv,z_2\inv,(z_1-z_2)\inv][[t]]$
such that \eqref{eq:expansionXisA} holds.
\end{proof}

Note that there need be no uniform $N$ that makes
\eqref{eq:rationality} true for all $k$; consider for instance the
case $A(z_1,z_2;t)=i_{z_1;z_2}e^{t/(z_1-z_2)}$.

\begin{lem} 
If there are 
\[
A(z_1, z_2;t)\in
 W((z_1))((z_2))[[t]],\quad B(z_2,z_1;t)\in
 W((z_2))((z_1))[[t]]
\]
such that there is for all $k\ge 0$ an $N\ge 0$ such that
   \begin{equation}
     \label{eq:generalizedlocality}
     (z_1-z_2)^NA(z_1z_2;t)\equiv
     (z_1-z_2)^NB(z_2,z_1;t)\mod t^k,
   \end{equation}
   then there is $X(z_1, z_2)\in
   W[[z_1,z_2]][z_1\inv,z_2\inv,(z_1-z_2)\inv][[t]]$ such that
 \[
 i_{z_1;z_2}X(z_1,z_2)=A(z_1,z_2),\quad i_{z_2;z_1}X(z_1,z_1)=B(z_2, z_1),
 \]
and
\[
(z_1-z_2)^N A(z_1,z_2)=(z_1-z_2)^N B(z_2,z_1)=(z_1-z_2)^N X(z_1,z_2).
\]
\end{lem}

\section{First Consequences of the Definition}
\label{sec:FrstConsq}

\begin{lem}\label{lem:Dvac=0}
  \[
D1=0.
\]
\end{lem}

\begin{proof}
  By the vacuum axiom \eqref{eq:vacuumX} for $i=1,2$ we have
\[
X_{z_1,z_2}(1\otimes 1)=e^{z_1D}1=e^{z_2D}1\in V[[z_1]]\cap V[[z_2]].
\]
This implies $D1=0$.
\end{proof}

We emphasize that $X_{z_1,z_2}$ is assumed to be nonsingular in the
$z_2$ variable at zero, so that $X_{z_1,0}$ is defined. (We used this
already in the locality axiom, (\ref{eq:localityAx}).) 
Define
\begin{equation}
  \label{eq:defXz}
  X_z\colon V\to V[[z]],\quad a\mapsto X_z(a)=e^{zD}a.
\end{equation}
We think of $X_z$ as the ``singular multiplication of $1$ element of
$V$'', which happens to be nonsingular, just as $X_{z_1,z_2}$ is the
singular multiplication of $2$ elements. Later, in Theorem
\ref{thm:analcontn}, we will define a singular multiplication
$X_{z_1,\dots,z_n}$ of $n$ elements of $V$.

Then we have
\begin{equation}
  \label{eq:relZztoXzz}
  X_z(a)=X_{z,0}(a\otimes 1),
\end{equation}
by the vacuum axiom \eqref{eq:vacuumX}.
\begin{lem}\label{lem:ExpX2inXX}
  For all $a, b\in V$ we have the following expansion:
\[
i_{z_1;z_2}X_{z_1,z_2}(a\ten b)=X_{z_1,0}(1\ten X_{z_2})(a\ten b).
\]
\end{lem}

\begin{proof}
  Since $X_{z_1,z_2}$ is regular at $z_2=0$ we have 
  \begin{align*}
i_{z_1;z_2}X_{z_1,z_2}(a\ten b)&=e^{z_2\partial_w}X_{z_1,w}\mid_{w=0}=\\
&=X_{z_1,0}(a\ten e^{z_2D}b) &&\text{ by \eqref{eq:H_DcovX},}=\\
&=X_{z_1,0}(1\ten X_{z_2})(a\ten b) &&\text{ by \eqref{eq:defXz}.}
\end{align*}
\end{proof}

\begin{remark}\label{rem:expansionversuscovariance}
  We derived the expansion of Lemma \ref{lem:ExpX2inXX} from the
  covariance axiom \eqref{eq:H_DcovX}. Conversely, if we know that
  $X_{z_1,z_2}$ has this expansion we see that
  $\partial_{z_2}X_{z_1,z_2}(a\ten b)$ and $X_{z_1,z_2}(a\ten Db)$
  both have the same image under $i_{z_1;z_2}$. So, $i_{z_1;z_2}$
  being injective, we can derive the covariance axiom
  \eqref{eq:H_DcovX} from the existence of the expansion in Lemma
  \ref{lem:ExpX2inXX}.
\end{remark}

\section{Analytic Continuation for $n=2$}
\label{sec:analn=2}
To make contact with the usual notation and terminology in the theory of vertex
algebras we introduce some definitions.

\begin{defn}[\textbf{Field}]\label{defn:field} Let $V$ be a
  $k$-module. A \emph{field} on $V$ is an element of $
  \Hom(V,V((z))[[t]])$.
\end{defn}

So if $a(z)$ is a field, we have for all $b\in V$ 
\[
a(z)b\in V((z))[[t]].
\]
\begin{defn}[\textbf{Vertex operator}]\label{def:vertexoperatorY}
  If $V$ is an $H_D$-quantum vertex algebra we define the vertex
  operator $Y(a,z)$ associated to $a\in V$ by
\begin{equation}
  \label{eq:DefY}
  Y(a,z)b=X_{z,0}(a\otimes b),
\end{equation}
  for $b\in V$. We will also use the notation $Y_z: a\otimes b\mapsto
  Y(a,z)b$, so that $Y_z=X_{z,0}$.
\end{defn}
Note that the vertex operator $a(z)=Y(a,z)$ for an $H_D$-quantum
vertex algebra is a field, for all $a\in V$.

We can rewrite Lemma \ref{lem:ExpX2inXX} as follows:
\begin{cor}[\textbf{Analytic continuation}]
\label{cor:AnalCont2}
The singular multiplication  $X_{z_1,z_2}(a\ten b)$ is the analytic continuation of
the product of vertex operators $Y(a,z_1)Y(b,z_2)1$, i.e.,
\[
i_{z_1;z_2}X_{z_1,z_2}(a\ten b)=Y(a,z_1)Y(b,z_2)1.
\]
\end{cor}
\begin{remark}
In Theorem \ref{thm:analcontn} we construct an $n$-variable version
$X_{z_1,z_2,\dots,z_n}$ of the singular multiplication satisfying
\[
i_{z_1;z_2;\dots;z_n} X_{z_1,z_2,\dots,z_n}(a_1\ten a_2\ten\dots\ten
a_n)=Y(a_1,z_1)Y(a_2,z_2)\dots Y(a_n,z_n)1,
\]
i.e., we construct the analytic continuation of arbitrary product of
vertex operators.
\end{remark}

\begin{remark}
  At this point we would like to emphasize that the axions we use are
  much weaker than those of Frenkel-Reshetikhin,
  \cite{q-alg/9706023}. Indeed, one of their axioms not only requires
  that the product of (quantum) vertex operators can be analytically
  continued, but also that the resulting function is meromorphic in the
  variables. This is not always the case in our $H_D$-quantum vertex
  algebras. For instance, we allow a singular multiplication
  $X_{z_1,z_2}(a\ten b)$ with a singularity of the form
  $e^{t/(z_1-z_2)}$, but this does not satisfy the Frenkel-Reshetikhin
  axioms, as there is an essential singularity at $z_1=z_2$. In our
  setup the quantum parameter $t$ is an independent variable (and we always
  expand in positive powers of $t$), whereas in Frenkel-Reshetikhin
  $t$ is a complex number.
\end{remark}
\section{Alternative Axioms}
\label{sec:AltAx}

We have formulated the axioms of an $H_D$-quantum vertex algebra in
terms of the rational singular multiplication $X_{z_1,z_2}$.
Traditionally the axioms of a vertex algebra have been formulated in
terms the 1-variable vertex operator $Y_z$. Let us briefly indicate how
this would work for $H_D$-quantum vertex algebras. Our axioms from
Definition \ref{defn:h_d-quantum-vertex-alg} would change slightly. We
start out with assuming the existence of a map
\[
Y_z\colon V\otimes V\to V((z))[[t]],
\]
instead of the singular multiplication $X_{z_1,z_2}$. The braiding and
translation maps $S^{(\tau)}$ and $S^{(\gamma)}$ are as before. The
vertex operator satisfies the following axioms:
\begin{description}
\item[(vacuum)] 
\[
Y_z(1\otimes a)=a,\quad Y_z(a\otimes 1)=e^{zD}a
\]
\item[($H_D$-covariance)]
\begin{align}
  \label{eq:CovarianceY2}i_{z;\gamma}Y(a,z+\gamma)e^{\gamma
    D}b&=i_{z;\gamma}e^{\gamma D}Y_{z}\circ S^{(\gamma)}_{z,0}(a\ten
  b).
\end{align}
\item[(Compatibility with Multiplication)] 
  \begin{align}
    S_{z_1,z_2}(Y_{w}\otimes 1)
    &=(Y_{w}\otimes1)i_{z_1, z_1-z_2;w}
    S_{z_1+w,z_2}^{23}S_{z_1,z_2}^{13},\label{eq:CompatYx1}\\
    S_{z_1,z_2}(1\otimes Y_{w})
    &=(1\otimes Y_{w})i_{z_1-z_2, z_2;w}
    S_{z_1,z_2+w}^{12}S_{z_1,z_2}^{13}.\label{eq:Compat1xY}
  \end{align}
\item[ (locality)] For all $a, b \in V$ and $k\ge 0$ there
exist $N$ such that for all $c\in V$
\begin{equation}
  \label{eq:LocalityY}
(z-w)^N Y(a,z)Y(b,w)c\equiv(z-w)^N Y_w(1\ten Y_z)(S_{w,z}(b\ten
a)\ten c) \mod t^k.
\end{equation}
\end{description}

Given these axioms we can reconstruct $X_{z_1,z_2}$.

\begin{lem}
  There exists a map
\[
X_{z_1,z_2}\colon V\ten V\to V[[z_1,z_2]][z_1\inv, (z_1-z_2)\inv][[t]]
\]
such that
\[
i_{z_1;z_2}X_{z_1,z_2}(a\ten b)= Y(a,z_1)Y(b,z_2)1.
\]
\end{lem}
\begin{proof}
Let $A(z_1,z_2)=  Y(a,z_1)Y(b,z_2)1$. By definition of the braiding
$S^{(\tau)}$ and the locality \eqref{eq:LocalityY} we have for all
$k\ge 0$ an $N\ge0$ such that
\[
(z_1-z_2)^NA(z_1,z_2)\in V[[z_1,z_2]][z_1\inv]][[t]] \mod t^k,
\]
and the Lemma follows from Lemma \ref{lem:z1-z2Arational}.
\end{proof}

Thus we can \emph{define}  in the present setup the singular
multiplication $X_{z_1,z_2}$ to be the analytic continuation of the product
$Y(a,z_1)Y(b,z_2)1$.

Alternatively, given the fields $Y(a,z)$ for any $a\in V$ we can
define $X_{z_1,z_2}$ as follows:
\begin{defn} \label{defn:XinY}For any $a, b\in V$ define
\[
X_{z_1,z_2}(a\ten b)=e^{z_2D}Y_{z_1-z_2}S^{(z_2)}_{z_1-z_2,0}(a\ten b).
\]
\end{defn}
The two definitions are equivalent:
\begin{lem} 
  \label{lem:ExpansionXfromYaxioms}If $X_{z_1,z_2}$ is given by
  Definition~\ref{defn:XinY} then
\[
   i_{z_1;z_2} X_{z_1,z_2}(a\ten b)=Y(a,z_1)e^{z_2D}b=Y(a,z_1)Y(b,z_2)1.
  \]
\end{lem}
The proof follows from (\ref{eq:CovarianceY2}) and the vacuum axiom.

To obtain the axioms of Section \ref{sec:H_DQuantumVertex} note that
Lemma \ref{lem:ExpansionXfromYaxioms} implies the covariance axiom
\eqref{eq:H_DcovX}, by Lemma \ref{lem:ExpX2inXX} and Remark
\ref{rem:expansionversuscovariance}. 
The rest of the axioms follow
immediately.

\begin{remark}
  We give in Definition \ref{defn:XinY} a direct construction of
  $X_{z_1,z_2}$ in terms of $Y_z$, without using analytic
  continuation. It seems not so easy to give such an explicit formula
  for the singular multiplication $X\subz$ of $n$ elements of $V$, to
  be introduced in Theorem \ref{thm:analcontn} using analytic
  continuation.
\end{remark}

\begin{remark}
  In the case of classical vertex algebras, as well as Etingof-Kazhdan
  (EK) quantum vertex operator algebras or Frenkel-Reshetikhin deformed chiral
  algebras, the translation map $S^{(\gamma)}_{z_1,z_2}$ is the
  identity, so that in this case
\begin{equation}
\label{eq:forskewsym2}
X_{z_1,z_2}(a\ten b)=e^{z_2D}Y(a,z_1-z_2)b\in
V[[z_1,z_2]][(z_1-z_2)\inv].
\end{equation}
In particular in these cases we can let $z_1=0$ as $X_{z_1,z_2}(a\ten
b)$ is not singular for $z_1=0$. That is no longer the case for the
examples of vertex operators connected to symmetric polynomials.
Therefore we have allowed for singular multiplication maps which are
singular in $z_1$ (but not in $z_2$, if we want to be able to define
$Y_z$ fields as above). It is possible to modify the theory further to
allow for singularities in both the variables, but we haven't yet
encountered examples which would call for such generalization.\qed
\end{remark}

The conclusion of this section is that we can start either with $Y_z$
or with $X_{z_1,z_2}$ as fundamental ingredient in the theory. Since
there are by now hundreds of papers on vertex algebras written in
terms of $Y_z$ we have allowed ourselves to emphasize $X_{z_1,z_2}$ in
this paper.

\section{Braiding and Skewsymmetry}
\label{sec:braidingskew}
                                                                                                 
An important fact of the theory of classical vertex algebras is that
the singular multiplication maps $X_{z_1,z_2}$ are ``commutative'',
i.e., we have for any $a, b\in V$
\begin{equation*}
 X_{z_1,z_2}(a\ten b)=X_{z_2,z_1}(b\ten a).
\end{equation*}
In the case of $H_D$-quantum vertex algebras the singular multiplication maps $X_{z_1,z_2}$ on
$V^{\ten 2}$ are no longer ``commutative'', but rather ``braided
commutative'', as shown by the next Lemma.

\begin{lem}[\textbf{Braided Symmetry}]\label{lem:X2braiding}  For any $a, b\in V$
\[
  X_{z_1,z_2}(a\ten b)=X_{z_2,z_1}S_{z_2,z_1}^{(\tau)}(b\ten a).
\]
\end{lem}
\begin{proof}
  Let $E=X_{z_1,z_2}(a\ten b)$,
  $F=X_{z_2,z_1}S_{z_2,z_1}^{(\tau)}(b\ten a)$.  We have
  \begin{equation}\label{eq:2}
\begin{aligned}
    i_{z_1;z_2}E&= X_{z_1,0}(1\ten X_{z_2})(a\ten
    b)=&&\text{by Lemma \ref{lem:ExpX2inXX}}\\
&=X_{z_1,0}(1\ten X_{z_2,0})(a\ten b\ten 1)&&\text{by \eqref{eq:defXz}}.
  \end{aligned}
\end{equation}
On the other hand, by the same calculation,
  \begin{equation}\label{eq:6}
    i_{z_2;z_1}F=X_{z_2,0}(1\ten X_{z_1,0})  (i_{z_2;z_1}S_{z_2,z_1}^{(\tau)}(b\ten a)\ten 1).
\end{equation}
By the locality axiom \eqref{eq:localityAx} the RHSs of
\eqref{eq:2} and \eqref{eq:6} are after multiplication by $(z_1-z_2)^N$
equal modulo $t^k$. But then there is for all $k\ge 0$ an $N\ge 0$
such that for the LHSs we have
\[
  (z_1-z_2)^N E \equiv (z_1-z_2)^N F \mod t^k.
\]
Since $E$ and $F$ both belong to $V[[z_1,z_2]][z_1\inv,z_2\inv
,(z_1-z_2)\inv][[t]]$ this implies that they are in fact equal.
\end{proof}
\begin{cor}[\textbf{Skewsymmetry}] \label{cor:skewsym}For any $a, b\in V$ we have
\[
e^{z_2D}Y_{z_1-z_2}\circ S_{z_1-z_2, 0}^{(z_2)}(a\ten b)
=e^{z_1D}Y_{z_2-z_1}\circ S_{z_2-z_1, 0}^{(z_1)}\circ
S_{z_2,z_1}^{(\tau)}(b\ten a).
\]
\end{cor}
The proof follows from Lemma \ref{lem:X2braiding} and Definition
\ref{defn:XinY}. 
\begin{remark} In the case of EK quantum vertex operator algebras the
  translation map $S^{(\gamma)}_{z_1,z_2}$ is the identity, \emph{and}
  the braiding map depends on a single variable $z_1-z_2$, therefore
  we can substitute $z_1=0$ and we get the EK braided skewsymmetry
  relation
\[
e^{zD}Y(a, -z)b= Y_{z}\circ S_{z}^{(\tau)}(b\ten a),
\]
where $S^{(\tau)}_z=S^{(\tau)}_{z,0}$.
 
Note that we cannot substitute $z_1=0$ in general as
$S^{(\tau)}_{z_2,z_1}$ might be singular at $z_1=0$, see Section
\ref{sec:mainexample}. The skewsymmetry relation in Corollary
\ref{cor:skewsym} looks much less appealing than the braided symmetry
relation in Lemma \ref{lem:X2braiding}. Many of the properties of
$H_D$-quantum vertex algebras look more symmetric in terms of the
singular maps $X_{z_1, z_2}$, which was one of the reasons we prefer
working with them, rather than the $Y_z$ fields.
\end{remark}

\section{Braiding Maps for $n>2$}
\label{sec:braidingmaps}

The singular multiplication map $X_{z_1,z_2}$ on $V^{\ten 2}$ is
invariant under simultaneous interchange of the variables $z_1,z_2$
and the factors in $V^{\ten 2}$, up to insertion of the two variable
braiding map $S_{z_1,z_2}^{(\tau)}$, according to the Lemma
\ref{lem:X2braiding}.  In the next section we will construct for all
$n\ge1$ a singular multiplication map $X\subz$ on $V^{\ten n}$, see
Theorem \ref{thm:analcontn}. These are invariant under simultaneous
permutation of the variables $z_i$ and the factors in $V^{\ten n}$, up
to insertion of an $n$ variable braiding map $S\subz^{f}$, see
Corollary \ref{cor:CompleteExpansionXinY}. In this section we
construct these braiding maps.

Let $n\ge 2$, $I_n=\{1,2,\dots,n \}$ and let $\mathcal{S}_n$ be
permutation group of $I_n$, i.e., the group of bijections
$\mathbf{f}\colon I_n\to I_n$.  Let $\w_i=(ii+1)\in \mathcal{S}_n$
(where $i=1,2,\dots,n$) be the simple transposition given on $j\in
I_n$ by
\[
\w_i(j)=
\begin{cases}
  j    & j\ne i,i+1\\
  i+1  & j=i\\
  i    & j=i+1
\end{cases}
\]
Then $\mathcal{S}_n$ is generated by the $\w_i$, with as only relations
\begin{equation}
  \label{eq:relSn}
\w_i^2=1,\quad \w_i\w_{i+1}\w_i=\w_{i+1}\w_i\w_{i+1},
\end{equation}
and
\begin{equation}
  \label{eq:relSntwo}
\w_i\w_j=\w_j\w_i,\quad \abs{i-j}\ge 2.
\end{equation}
Now let $V$ be a free $k$-module, and define a \emph{right} action for
$\mathbf{f}\in \mathcal{S}_n$ on the $n$-fold tensor product of $V$ by
\[
\sigma_\mathbf{f}\colon V^{\ten n}\to V^{\ten n},\quad A_n\mapsto
a_{\mathbf{f}(1)}\ten a_{\mathbf{f}(2)}\ten \dots a_{\mathbf{f}(n)},
\]
where $A_n=a_1\ten a_2\ten \dots a_n\in V^{\ten n}$. Let $\tau\colon
a\ten b \mapsto b\ten a\in V^{\ten 2}$. Then 
\[
\sigma_{\w_i}=i^{i-1}\ten \tau\ten 1^{n-i-1}.
\]
Here we write $1^k$ for $1_V\ten1_V\ten\dots\ten1_V$, the $k$-fold
tensor product of the identity $1_V\colon V\to V$. 
We emphasize that if $\mathbf{f}=\mathbf{g}\w_i$ then
$\sigma_{\mathbf{f}}=\sigma_{\w_i}\sigma_{\mathbf{g}}$.

Let $\Ratz$ be a space of rational functions in $n$ variables. Then
$\mathcal{S}_n$ acts on the \emph{left} on $\Ratz$ by permutation the
variables: if $\mathbf{f}\in\mathcal{S}_n$ and $A\subz\in \Ratz$, then
we put
\[
\mathbf{f}.A\subz=A_{\mathbf{f}(z_1,z_2,\dots,z_n)},
\]
where we write $\mathbf{f}(z_1,z_2,\dots,z_n)$ for
$z_{\mathbf{f}(1)},z_{\mathbf{f}(2)},\dots,z_{\mathbf{f}(n)}$. 

Now let $\MapzV$ be the space of linear maps
\[
V^{\ten n}\to V^{\ten n}[z_i^{\pm 1},(z_i-z_j)\inv][[t]],\quad 1\le
i<j\le n.
\]
We have an action of $\mathcal{S}_n$ on $\MapzV$ combining the action
of $\mathcal{S}_n$ on $V^{\ten n}$ and on rational functions:
if $\mathbf{f}\in \mathcal{S}_n$ and
$A\subz\in\MapzV$ then define
\[
\mathbf{f}.A\subz=\sigma_\mathbf{f}\inv\circ A_{\mathbf{f}(z_1,z_2,\dots,z_n)}\circ
\sigma_\mathbf{f}.
\]
Now let $V$ be an $H_D$-quantum vertex algebra. So we get, by definition, in
particular a braiding map $S^{(\tau)}_{z_1,z_2}\in
\Map_{z_1,z_2}(V^{\ten 2})$. For simplicity we denote it by
$S_{z_1,z_2}$ in this section, as we will not use
$S_{z_1,z_2}^{(\gamma)}$ here. It satisfies, see Definition
\ref{defn:h_d-quantum-vertex-alg},
\begin{align}\label{eq:propertiesS1}
S_{z_1,z_2}\circ\tau\circ S_{z_2,z_1}\circ \tau&=1_{V^{\ten 2}},\\
S_{z_1,z_2}^{12}S_{z_1,z_3}^{13}S_{z_2,z_3}^{23}&=S_{z_1,z_3}^{14}S_{z_1,z_3}^{13}S_{z_1,z_2}^{12}.
\label{eq:propertiesS2}
\end{align}
We will to use the braiding matrix $S_{z_1,z_2}$ to define a map
$\mathcal{S}_n\to \MapzV$.

\begin{defn}[\textbf{Braiding maps}]\label{defn:braiding-maps-n2}
  Define for each $\mathbf{f}\in\mathcal{S}_n$ an element $S\subz^{f}$
  of $\MapzV$,
  callled the braiding map associated to $\mathbf{f}$, by expanding
  $\mathbf{f}$ (in any way) in simple reflections $w_i$ and using
\[
S^{\w_i}\subz=1^{i-1}\ten S^{(\tau)}_{z_i,z_{i+1}}\ten i^{n-i-1},
\]
and
\begin{equation}
\label{eq:Sfg}
S^{\mathbf{f}{\mathbf{g}}}\subz=S^{\mathbf{g}}\subz \sigma_{\mathbf{g}}
S^\mathbf{f}_{{\mathbf{g}}\inv(z_1,\dots,z_n)}(\sigma_\mathbf{f})^{-1}.
\end{equation}
\end{defn}
                                                                                                 
The point is that to define $S\subz^{\mathbf{f}}$ for
$\mathbf{f}\in\mathcal{S}_n$, we can take \emph{any} decomposition of
$\mathbf{f}$ into simple transpositions $\w_i$, i.e., this definition
is unambigious. The proof of this statement is discussed in Appendix
\ref{app:braiding}.

\section{Analytic Continuation for $n>2$}
\label{sec:Analcont}

If $V$ is an $H_D$-quantum vertex algebra, recall that we have the
``singular'' multiplications $X_z$ and $X_{z_1,z_2}$ of $1$,
respectively $2$ elements of $V$, see (\ref{eq:defXz}) and Definition
\ref{defn:h_d-quantum-vertex-alg}. We will in this section construct
singular multiplications $X\subz$ of $n$ elements of $V$.

Let ${\mathbf{f}}_n=\w_1\w_2\dots \w_{n-1}$ be the $n$-cycle $(123\dots n)$
and consider the associated braiding matrix $S^{{\mathbf{f}}_n}\subz$. We have
${\mathbf{f}}_n=\w_1(23\dots n)$. Writing ${\mathbf{f}}_{n-1}=(23\dots n)$ and
$\sigma_n=\sigma_{{\mathbf{f}}_n}$ we find from \eqref{eq:Sfg} that
\begin{equation}
  \label{eq:DefSsigman}
  S^{{\mathbf{f}}_n}_{z_2,z_3\dots,z_n,z_1}\sigma_{n}=(1\ten
  S^{{\mathbf{f}}_{n-1}}_{z_3,\dots,z_n,z_1}\sigma_{n-1})(S_{z_2,z_1}^{(\tau)}\tau
  \ten 1^{\ten n-2}).
\end{equation}
We will frequently use the abbreviation
\begin{equation}
  \label{eq:defp_n}
p_n=p_n(z_1,z_2,\dots,z_n)=\prod_{1\le i<j\le n} z_i-z_j.
\end{equation}

\begin{thm}[\textbf{Analytic Continuation}]\label{thm:analcontn} Let $V$ be an $H_D$-quantum vertex
  algebra. There exists for all $n\ge 2$ maps
\[
X_{z_1,\dots,z_n}\colon V^{\ten n}\to
V[[z_k]][z_j\inv,(z_i-z_j)\inv][[t]],\quad {1\le i<j \le n}, i\le k\le n
\]
such that
\begin{equation}
i_{z_1;z_2,\dots,z_n}X_{z_1,\dots,z_n}=X_{z_1,0}(1\ten
X_{z_2,\dots,z_n}).\label{eq:X1tenXexpansionXn}
\end{equation}
and
\begin{equation}
  X_{z_1,z_2,\dots,z_n}=X_{z_2,\dots,z_n,z_1}S^{{\mathbf{f}}_n}_{z_2,\dots,z_n,z_1}  \sigma_n,\label{eq:cyclicsymX}
\end{equation}
where $S^{{\mathbf{f}}_n}_{z_2,\dots,z_n,z_1}$ is defined in
Definition \ref{defn:braiding-maps-n2}.
\end{thm}
\begin{proof}
The theorem is true for $n=2$ by Lemma \ref{lem:ExpX2inXX} and Lemma
\ref{lem:X2braiding}. Assume the theorem is true for $\ell$, $2\le
\ell\le n_0$ and let $n=n_0+1$. The induction hypothesis implies that
\[
i_{z_2;z_3,\dots,z_n}X_{z_2,z_3,\dots,z_n}=X_{z_2,0}(1\ten
X_{z_3,z_4,\dots,z_n}),
\]
so that for all $k\ge0$ there is an $N\ge0$ such that
\begin{equation}
\label{eq:7}
p_{n-1}^N X_{z_2,z_3,\dots,z_n}\equiv p_{n-1}^NX_{z_2,0}(1\ten
X_{z_3,z_4,\dots,z_n})\mod t^k.
\end{equation}
Also we have 
\begin{equation}
  \label{eq:8}
  X_{z_1,z_3,z_4\dots,z_n}=X_{z_3,z_4,\dots,z_n,z_1}S^{{\mathbf{f}}_{n-1}}_{z_3,z_4,\dots,z_n,z_1}
  \sigma_{n-1}.
\end{equation}
Consider $E=X_{z_1,0}(1\ten X_{z_2,z_3,\dots,z_n})(A_n)$,
$A_n\in V^{\ten n}$. This is an
element of\\$V((z_1))[[z_2,z_3,\dots,z_{n}]][z_1\inv,\dots,z_{n-1}\inv,(z_i-z_j)\inv][[t]]$,
$2\le i<j\le n$, and we want to show $E$ is in the image of the
expansion $i_{z_1;z_2,\dots,z_n}$. For this it suffices to show that
for all $k\ge0$ there is an $N\ge0$ such that $p_n^N E$ has at worst a
finite order pole in $z_1$. This is a small calculation: for all $k\ge
0$ there is $N\ge0$ such that modulo $t^k$ we have
\begin{align*}
  &p_n^N E=p_n^N X_{z_1,0}(1\ten
  (X_{z_2,0}(1\ten X_{z_3,\dots,z_n})))(A_n)= &&\text{ by \eqref{eq:7}}\\
  &= p_n^N X_{z_2,0}(1\ten (X_{z_1,0}(1\ten
  X_{z_3,\dots,z_n})))(i_{z_2;z_1}S^{(\tau)}_{z_2,z_1}\tau\ten 1^{
    n-2})(A_n)= &&\text{ by \eqref{eq:localityAx}}\\
  &= p_n^N X_{z_2,0}(1\ten
  (X_{z_1,z_3,\dots,z_n})))(S^{(\tau)}_{z_2,z_1}\tau\ten 1^{n-2})(A_n)= &&\text{ by
    \eqref{eq:7} }\\
  &= p_n^N X_{z_2,0}(1\ten (X_{z_3,z_4,\dots,z_n,z_1})))(1\ten
  S^{{\mathbf{f}}_{n-1}}_{z_3,z_4,\dots,z_n,z_1}\sigma_{n-1})\times\\
  &\qquad\qquad \times (S^{(\tau)}_{z_2,z_1}\tau\ten 1^{n-2})(A_n)= &&\text{
    by \eqref{eq:8} }\\
  &= p_n^N X_{z_2,0}(1\ten (X_{z_3,z_4,\dots,z_n,z_1})))
  S^{{\mathbf{f}}_n}_{z_2,z_3,\dots,z_n,z_1}\sigma_n(A_n)= &&\text{ by
    \eqref{eq:DefSsigman}}.
\end{align*}
We see from the last expression that $p_n^N E$ has at worst a finite
order pole in $z_1$ and hence there is $X_{z_1,z_2,\dots,z_n}$ such
that \eqref{eq:X1tenXexpansionXn} holds.

Next consider
$F=X_{z_2,z_3,\dots,z_n,z_1}S^{f_n}_{z_2,z_3,\dots,z_n,z_1}\sigma_n
(A_n)$ and $G=X\subz(A_n)$. For all $k\ge 0$ there is an $N\ge0$ such
that modulo $t^k$
\[
p_n^N E=p_n^N G.
\]
By what we just proved we have
\[
i_{z_2;z_3,\dots,z_n,z_1}F=X_{z_2,0}(1\ten X_{z_3,\dots,z_n,z_1})i_{z_2;z_3,\dots,z_n,z_1}
S^{{\mathbf{f}}_n}_{z_2,z_3,\dots,z_n,z_1}\sigma_n(A_n),
\]
and so for all $k\ge
0$ there is $N\ge0$ such that modulo $t^k$ we have
\[
p_n^N F=p_n^N X_{z_2,0}(1\ten
X_{z_3,\dots,z_n,z_1})i_{z_2;z_3,\dots,z_n,z_1}S^{{\mathbf{f}}_n}_{z_2,z_3,\dots,z_n,z_1}\sigma_n
(A_n)=p_n^N E=p_n^N G.
\]
Since $G$, $F$ both belong to $V[[z_i]][z_i\inv, (z_i-z_j)\inv][[t]]$
this forces $G=F$, i.e., \eqref{eq:cyclicsymX} holds. 
\end{proof}

\begin{cor}[\textbf{Analytic Continuation for products of fields}]\label{cor:CompleteExpansionXinY}
  For all $n\ge2$ and $1\le i\le n-1$ we have, if $A_n=a_1\ten a_2\ten
  \dots \ten a_n$, the expansion
  \begin{multline*}
    i_{z_1;z_2;\dots,z_i;z_{i+1},z_{i+2},\dots,z_n}X_{z_1,\dots,z_n}(A_n)=X_{z_1,0}(a_1\ten
    X_{z_2,0}(a_2\ten\dots \\\dots (a_{i-1}\ten X_{z_i,0}(a_i\ten
    X_{z_{i+1},z_{i+2},\dots,z_n}(a_{i+1}\ten a_{i+2}\ten \dots\ten
    a_n)\dots)))).
\end{multline*}
\end{cor}

In particular the case of $i=n-1$ of the Corollary is (using the
notation \eqref{eq:DefY})
\begin{equation}
  \label{eq:FullExpansX}
i_{z_1;z_2;\dots;z_n}X_{z_1,z_2,\dots,z_n}(A_n)=Y(a_1,z_1)Y(a_2,z_2)\dots Y(a_n,z_n)1.
\end{equation}
In other words the $n$-variable $X_{z_1,z_2,\dots,
  z_n}(A_n)$ is the analytic continuation of the composition of $n$
vertex operators acting on the vacuum.

\section{Further Consequences}
\label{sec:furtcons}
\begin{lem}\label{lem:infH_Dcovariance}
  The infinitesimal forms of the $H_D$-covariance axioms
  \eqref{eq:H_DcovS}, \eqref{eq:H_DcovMult} are
\begin{align}
     &\label{eq:infScov}
( 1\ten D+\partial_{z_2})S_{z_1,z_2}=S_{z_1,z_2}(1\ten D),\\
     &  DY(a,z)b=\partial_z Y(a,z)b+Y(a,z)Db-Y_{z}\circ \alpha_{z,0}(a\otimes b).
 \end{align}
where $\alpha_{z_1,z_2}$ is defined to be $\partial_{\gamma }S^{(\gamma)}_{z_1,z_2}$.
and satisfies the infinitesimal form of the vacuum axioms \eqref{eq:vacuumS}
\[
\alpha_{z_1,z_2}(a\otimes 1)=0, \quad \alpha_{z_1,z_2}(1\otimes     b)=0,
\]
\end{lem}

\begin{lem}\label{lem:H_Dcovfirst}
For all $a,b\in V$
  \begin{equation}\label{eq:4}
    X_{z_1,z_2}(Da\otimes b)=\partial_{z_1}X_{z_1,z_2}(a\otimes b),
  \end{equation}
    \end{lem}

\begin{proof}
By the previous Lemma \ref{lem:infH_Dcovariance}
\begin{align*}
  &X_{z_1,z_2}(Da\ten b)=X_{z_1,z_2} \tau (1\ten D)(b\ten a)=&&
  \\
  &=X_{z_2,z_1}S^{(\tau)}_{z_2,z_1}(1\ten D)(b\ten a) =&&\text{by
    Lemma \ref{lem:X2braiding}}\\
  &=X_{z_2,z_1}(1\ten D+\partial_{z_2})S^{(\tau)}_{z_2,z_1}(b\ten a)
  =&&\text{by \eqref{eq:infScov}}\\
  &=\left(\partial_{z_2}(X_{z_2,z_1})S^{(\tau)}_{z_2,z_1}+X_{z_2,z_1}
  \partial_{z_2}S^{(\tau)}_{z_2,z_1}\right)(b\ten a) =&&\text{by
    \eqref{eq:H_DcovX}}\\
  &=\partial_{z_2}(X_{z_2,z_1}S^{(\tau)}_{z_2,z_1})(b\ten a) =\\
  &=\partial_{z_2}X_{z_1,z_2}(a\ten b)&&\text{by Lemma
    \ref{lem:X2braiding},}
\end{align*}
proving \eqref{eq:4}.
\end{proof}

\begin{cor}
  For all $n\ge 2$, $1\le i\le n$ and $A_n=a_1\ten a_2\ten\dots\ten
  a_n\in V^{\ten n}$ we have
\[
\partial_{z_i}X_{z_1,z_2,\dots,z_n}(A_n)=X_{z_1,z_2,\dots,z_n}(a_1\ten\dots\ten
Da_i\ten\dots\ten a_n)).
\]
\end{cor}
\begin{proof}
  For $n=2$ and $i=2$ this is axiom \eqref{eq:H_DcovX}, and for $i=1$
  this is Lemma \ref{lem:H_Dcovfirst}.  Put for $n>2$
\[
E=\partial_{z_i}X_{z_1,z_2,\dots,z_n}(A_n),\quad F=X_{z_1,z_2,\dots,z_n}(a_1\ten\dots\ten
Da_i\ten\dots\ten a_n).
\]
Since $\partial_{z_i}$ commutes with expansions we have
\begin{align*}
  & i_{z_1;z_2;\dots;z_i;z_{i+1},\dots,z_n}E=
  X_{z_1,0}(a_1\ten  X_{z_2,0} (a_2\ten\dots \\
  &\dots(a_{i-1}\ten\partial_{z_i}X_{z_i,0}(a_i\ten
  X_{z_{i+1},\dots,z_n}(a_{i+1}\ten\dots\ten a_n))\dots)))=&&\text{ by
    Corollary \ref{cor:CompleteExpansionXinY}}\\
  &=X_{z_1,0} (a_1\ten  X_{z_2,0}(a_2\ten\dots\\
  &\dots (a_{i-1}\ten X_{z_i,0}(Da_i\ten X_{z_{i+1},\dots,z_n}(a_{i+1}\ten\dots\ten
  a_n)))\dots))=&&\text{ by
    Lemma \ref{lem:H_Dcovfirst}}\\
  &=i_{z_1;z_2;\dots;z_i;z_{i+1},\dots,z_n}F&&\text{ by Corollary
    \ref{cor:CompleteExpansionXinY}}
\end{align*}
Since both $E$ and $F$ belong to $V[[z_i]][z_i\inv,(z_i-z)\inv][[t]]$,
$1\le i<j\le n$, and have the same expansion they must be equal.
\end{proof}

\begin{lem} \label{lem:VacnX}
  For all $n\ge 1$ and  $A_n=a_1\ten a_2\ten\dots\ten a_n\in V^n$ we have
   \[
   X_{z_1,z_2,\dots, z_n,0}(A_n \otimes   1)=X_{z_1,z_2,\dots,z_n}(A_n).
   \]
\end{lem}

\begin{proof}
  For $n=1$ this is \eqref{eq:relZztoXzz}. Assume that the lemma is
  true for all $\ell$, $1\le \ell\le n_0$, and let $n=n_0+1$. Put
  $E=X_{z_1,z_2,\dots,z_n,0}(A_n\ten 1)$,
  $F=X_{z_1,z_2,\dots,z_n}(A_n)$. By Theorem \ref{thm:analcontn} and
  the induction hypothesis
\begin{align*}
  i_{z_1;z_2,\dots,z_n}E&= X_{z_1,0}(a_1\ten
  X_{z_2,\dots,z_n,0}(a_2\ten a_3\ten\dots \ten a_n\ten 1))=\\
  &=X_{z_1,0}(a_1\ten X_{z_2,\dots,z_n}(a_2\ten a_3\ten\dots \ten
  a_n))=\\
  &=i_{z_1;z_2,\dots,z_n}F.
\end{align*}
Since both $E$ and $F$ belong to $V[[z_i]][z_i\inv,(z_i-z)\inv][[t]]$,
$1\le i<j\le n$, and have the same expansion they must be equal.
\end{proof}

\begin{lem}  Suppose that $S^{(\gamma)}_{z_1,z_2}$ is the
    identity map on $V\otimes V$. Then the
    following is true:
\begin{gather}
\label{eq:1}DX_{z_1, z_2}=(\partial_{z_1}+\partial_{z_2})X_{z_1, z_2}=X_{z_1, z_2}(D\ten 1+1\ten D),\\
\label{eq:9}[D,Y(a,z)]=\partial_z Y(a,z),\\
\label{eq:10}X_{z_1,z_2}\circ (\partial_{z_1}+\partial_{z_2}) S^{(\tau)}_{z_1,z_2}=0.
\end{gather}

\end{lem}
\begin{proof}
  The second property is a direct consequence of Lemma
  \ref{lem:infH_Dcovariance}. The first equality follows from
  expanding both sides of \eqref{eq:H_DcovMult} in powers of $\gamma$
  and comparing the coefficients in front of $\gamma ^1$. 

For the last part rewrite \eqref{eq:1} as 
\[
e^{\gamma D}X_{z_1,z_2}=X_{z_1,z_2}\Delta(e^{\gamma D}),
\]
where $\Delta$ is the coproduct of $H_D$, so that $\Delta(e^{\gamma
  D}=e^{\gamma D}\ten e^{\gamma D})$. 
Similarly rewrite the
$H_D$-covariance axiom \eqref{eq:H_DcovS} for the braiding as
\[
(1\ten e^{-\gamma D})
S_{z_1,z_2}=e^{\gamma(\partial_{z_1}+\partial_{z_2})}S_{z_1,z_2}(1\ten
 e^{-\gamma D}).
\]
By differentiating with respect to $z_1$ the axiom \eqref{eq:GrpSStau}
we obtain a similar equation involving $\partial_{z_1}$ and
$e^{-\gamma D}\ten 1$, and we combine these as
\[
\Delta(e^{-\gamma D}) S^{(\tau)}_{z_1,z_2}=
e^{\gamma(\partial_{z_1}+\partial_{z_2})}
S^{(\tau)}_{z_1,z_2}\Delta(e^{-\gamma D}) .
\]
Now we calculate
\begin{align*}
  e^{-\gamma D} X_{z_2,z_1}&=   e^{-\gamma D} X_{z_1,z_2}S_{z_1,z_2}\tau=\\
  &= X_{z_1,z_2}\Delta(e^{-\gamma D})S_{z_1,z_2}\tau=\\
  &= X_{z_1,z_2}e^{\gamma(\partial_{z_1}+\partial_{z_2})}
S^{(\tau)}_{z_1,z_2}\tau\Delta(e^{-\gamma D}).
\end{align*}
On the other hand
\begin{align*}
  e^{-\gamma D} X_{z_2,z_1}&=  X_{z_1,z_2}S_{z_1,z_2}\tau\Delta(
  e^{-\gamma D} ).
\end{align*}
By multiplying by $\Delta(e^{\gamma D})\tau$ on the right we find
\[
X_{z_1,z_2}e^{\gamma(\partial_{z_1}+\partial_{z_2})}
S^{(\tau)}_{z_1,z_2}=
X_{z_1,z_2}S^{(\tau)}_{z_1,z_2},
\]
from which \eqref{eq:10} follows.
\end{proof}
\begin{remark}
  In the context of the lemma above it is natural to assume that
  $S^{(\tau)}_{z_1,z_2}$ is a function of just $z_1-z_2$. In this case
  $V$ is a quantum vertex operator algebra as defined by
  Etingof-Kazhdan, see \cite{MR2002i:17022} (except for the fact that
  they insist that the braiding is of the form
  $S^{(\tau)}=1+\mathcal{O}(t)$).
\end{remark}

\section{Braiding and singular multiplication}
\label{sec:brSingMult}

We have seen that the $n$-fold singular multiplication has cyclic
symmetry: if ${\mathbf{f}}_n$ is the cyclic permuation $(123\dots n)$, then
\begin{equation}
  \label{eq:cyclicsymXn}
  X\subz=  X_{{\mathbf{f}}_n(z_1,\dots,z_n)}S^{{\mathbf{f}}_n}_{{\mathbf{f}}_n(z_,\dots,z_n)}
  \sigma_n,
\end{equation}
  see Theorem \ref{thm:analcontn}. In this section we show that in fact
the $n$-fold singular multiplication has arbitrary permutation
symmetry: in \eqref{eq:cyclicsymXn} we can replace ${\mathbf{f}}_n$ by any ${\mathbf{f}}\in \mathcal{S}_n$.

\begin{lem}\label{lem:transpositiontau}
 For all $n\ge 2$ we have
\[
X_{z_1,z_2,\dots,z_n}=X_{\w_1(z_1,z_2,\dots,z_n)} S^{\w_1}\subz
(\tau\ten 1^{n-2}).
\]  
\end{lem}
\begin{proof}Let $E=X_{z_1,z_2,\dots,z_n}(A_n)$,
  $F=X_{z_2,z_1,\dots,z_n}\circ (S^{(\tau)}_{z_2,z_1}\tau\ten 1^{n-2}) (A_n)$,
  $A_n\in V^{\ten n}$. Then
  there exist for all $k\ge 0$ an $N\ge0$ such that modulo $t^k$
 \begin{align*}
   p_n^N&i_{z_1;z_2;z_3,z_4\dots,z_n}E=p_n^N X_{z_1,0}(1\ten
   X_{z_2,0}(1\ten X_{z_3,z_4,\dots,z_n}))(A_n)&&\text{by Thm \ref{thm:analcontn}}\\
&=p_n^N X_{z_2,0} (1\ten
   X_{z_1,0}(1\ten X_{z_3,z_4,\dots,z_n})) (S^{(\tau)}_{z_2,z_1}\ten
   1^{\ten n-2})(A_n)=&&\text{by \eqref{eq:localityAx}}\\
&=p_n^Ni_{z_1;z_2;z_3,z_4\dots,z_n}F,
 \end{align*}
 by Theorem \ref{thm:analcontn} again.  Since both $E$ and $F$ belong
 to $V[[z_i]][z_i\inv,(z_i-z)\inv][[t]]$, $1\le i<j\le n$, and have
 the same expansion they must be equal.
\end{proof}

Recall that the first simple transposition $\w_1$ and the cyclic
permutation $f_n=(123\dots n)$ generate $\mathcal{S}_n$.

\begin{cor}\label{cor:permsymmetryX}
  If ${\mathbf{f}}\in \mathcal{S}_n$ is a permutation of $\{1, 2,\dots,n\}$ and
  $\sigma_{\mathbf{f}}(A_n)=a_{{\mathbf{f}}(1)}\ten a_{{\mathbf{f}}(2)}\ten \dots\ten a_{{\mathbf{f}}(n)}$, then
  \begin{equation}
    \label{eq:Xfsym}
    X\subz=X_{{\mathbf{f}}(z_{1},\dots,z_{n})}
    S^{{\mathbf{f}}}_{{\mathbf{f}}(z_{1},
      \dots,z_{n})}\sigma_{\mathbf{f}}.
  \end{equation}
   \end{cor}

\begin{proof}
  Suppose we have two elements ${\mathbf{f}},\mathbf{g}\in \mathcal{S}_n$ such that
  \eqref{eq:Xfsym} holds. Then, by \eqref{eq:multphifg},
  \begin{align*}
    X_{{\mathbf{f}}\mathbf{g}(z_1,\dots,z_n)}S^{{\mathbf{f}}\mathbf{g}}_{{\mathbf{f}}\mathbf{g}(z_1,\dots,z_n)}\sigma_{{\mathbf{f}}\mathbf{g}}&=
    X_{{\mathbf{f}}\mathbf{g}(z_1,\dots,z_n)}S^{\mathbf{g}}_{{\mathbf{f}}\mathbf{g}(z_1,\dots,z_n)}\sigma_{\mathbf{g}}
    S^{\mathbf{f}}_{\mathbf{g}(z_1,\dots,z_n)}\sigma_{{\mathbf{f}}}=\\
&={\mathbf{f}}.(X_{\mathbf{g}(z_1,\dots,z_n)}S^{\mathbf{g}}_{\mathbf{g}(z_1,\dots,z_n)}\sigma_{\mathbf{g}})
    S^{\mathbf{f}}_{{\mathbf{f}}(z_1,\dots,z_n)}\sigma_{{\mathbf{f}}}=\\
&={\mathbf{f}}.(X_{(z_1,\dots,z_n)})S^{\mathbf{f}}_{{\mathbf{f}}(z_1,\dots,z_n)}\sigma_{{\mathbf{f}}}=\\
&=X\subz.
  \end{align*}
So if \eqref{eq:Xfsym} holds for ${\mathbf{f}}$ and for $\mathbf{g}$ it holds for
${\mathbf{f}}\mathbf{g}$. But we know that \eqref{eq:Xfsym} holds for ${\mathbf{f}}_n$, by Theorem
\ref{thm:analcontn}, and for $\w_1$, by Lemma
\ref{lem:transpositiontau}, and these elements generate
$\mathcal{S}_n$. So \eqref{eq:Xfsym} holds for all ${\mathbf{f}}\in \mathcal{S}_n$.
\end{proof}
\section{Expansions of $X_{z_1,z_2,0}$}
\label{sec:ExpX}

We have seen that the expansion of $X\subz$ (in the region $\abs
{z_1}>\abs{z_2}>\dots>\abs{z_n}$) is expressed as a
composition of $1$-variable vertex operators. In particular, for $n=3$
we get, if $A=a\ten b\ten c$,
\[
i_{z_1;z_2}X_{z_1,z_2,0}(A)=Y(a,z_1)Y(b,z_2)c,
\]
see \eqref{eq:FullExpansX}. In this section we find other expansions
of $X_{z_1,z_2,0}$ that have useful expressions in terms of $Y_z$.

First we need a variant of the analytic continuation Theorem
\ref{thm:analcontn}. 

\begin{lem}
  \[
X_{z_1,z_2}(1\ten X_{w,0}i_{z_2;w}S^{(z_2)}_{w,0})=i_{z_1-z_2,z_2;w}X_{z_1,z_2+w,z_2}.
\]
\end{lem}

\begin{proof}
  \begin{align*}
i_{z_1;z_2}& X_{z_1,z_2}(1\ten X_{w,0}i_{z_2;w}S^{(z_2)}_{w,0})=\\
&=X_{z_1,0}(1\ten e^{z_2D}X_{w,0}i_{z_2;w}S^{(z_2)}_{w,0})=&&\text{ by
Lemma \ref{lem:ExpX2inXX}}\\
&=X_{z_1,0}(1\ten i_{z_2;w}X_{w+z_2,z_2})=&&\text{ by Axiom \eqref{eq:H_DcovMult}}\\
&=i_{z_2;w}i_{z_1;w+z_2,z_2}X_{z_1,z_2+w,z_2}=&&\text{ by Thm. \ref{thm:analcontn}}\\
&=i_{z_1;z_2}i_{z_1-z_2,z_2;w}X_{z_1,z_2+w,z_2},
  \end{align*}
since
\[
i_{z_1;z_2}i_{z_1-z_2;w}f(z_1-z_2-w)=i_{z_2;w}i_{z_1;w+z_2}f(z_1-z_2-w).
\]
The Lemma follows then by cancelling $i_{z_1;z_2}$.
\end{proof}

Next we need a variant of the compatibility with multiplication Axiom
\eqref{eq:Compat1xX} .

\begin{lem}\label{lem:compatwithXS}

\[
  S^{(\tau)}_{z_1,z_2} (1\ten X_{w,0}S_{w,0}^{(\gamma)})=(1\ten
  X_{w,0}S_{w,0}^{(\gamma)}) i_{z_1-z_2,z_2;w}S^{(\tau)12}_{z_1;z_2+w}
  S^{(\tau)13}_{z_1,z_2}.
\] 
\end{lem}

\begin{proof}
  We need some simple identities. By Axiom \eqref{eq:H_DcovS}
  \begin{equation}
    \label{eq:a}
X_{w,0}S^{(\gamma)}_{w,0}=e^{-\gamma D}X_{w+\gamma,\gamma}.
  \end{equation}
By Axiom \eqref{eq:H_DcovS}
\begin{equation}
  \label{eq:b}
S^{(\tau)}_{z_1,z_2}(1\ten e^{-\gamma D})=(1\ten e^{-\gamma
  D})i_{z_1-z_2,z_2;\gamma}S^{(\tau)}_{z_1,z_2-\gamma}.
\end{equation}
Finally, by Axiom \eqref{eq:Compat1xX}
\begin{equation}
  \label{eq:c}
S^{(\tau)}_{z_1,z_2}(1\ten X_{w,0})=(1\ten X_{w,0})i_{z_1-z_2,z_2;w}
S^{(\tau)12}_{z_1,z_2+w}S^{(\tau)13}_{z_1,z_2}.
\end{equation}
Then
\begin{align*}
  S^{(\tau)}_{z_1,z_2}&(1\ten
  X_{w,0}S_{w,0}^{(\gamma)})=S^{(\tau)}_{z_1,z_2}(1\ten e^{-\gamma{
      D}}X_{w+\gamma,\gamma})=&&\text{by \eqref{eq:a}}\\
  &=(1\ten e^{-\gamma D}) i_{z_1-z_2,z_2;\gamma}
  S^{(\tau)}_{z_1,z_2-\gamma}
  (1\ten   X_{w+\gamma,\gamma})=&&\text{by \eqref{eq:b}}\\
  &=(1\ten e^{-\gamma D}) (1\ten   X_{w+\gamma,\gamma})\times\\
&\quad \times i_{z_1-z_2,z_2;\gamma}i_{z_1-z_2-\gamma,z_2-\gamma;w+\gamma,\gamma}
  S^{(\tau)12}_{z_1,z_2-\gamma+(w-\gamma)}  S^{(\tau)13}_{z_1,(z_2-\gamma)+\gamma}
  =&&\text{by \eqref{eq:c}}\\
&=(1\ten X_{w,0}S^{(\gamma)}_{w,0})\times\\
&\quad \times i_{z_1-z_2,z_2;\gamma}i_{z_1-z_2-\gamma,z_2-\gamma;w+\gamma,\gamma}
  S^{(\tau)12}_{z_1,z_2-\gamma+(w-\gamma)}  S^{(\tau)13}_{z_1,(z_2-\gamma)+\gamma}
  =&&\text{by \eqref{eq:a}}\\
&=(1\ten X_{w,0}S^{(\gamma)}_{w,0}) i_{z_1-z_2,z_2;w}  S^{(\tau)12}_{z_1,z_2+w}  S^{(\tau)13}_{z_1,z_2},
  \end{align*}
since 
\begin{align*}
i_{z_1-z_2-\gamma,z_2-\gamma;w+\gamma,\gamma}f((z-\gamma)+(w+\gamma))&=i_{z_1-z_2,z_2;w}f(z+w)\\
i_{z_1-z_2-\gamma,z_2-\gamma;\gamma}f((z-\gamma)+\gamma)=f(z).
\end{align*}
\end{proof}

\begin{remark}
  Note that in Lemma \ref{lem:compatwithXS} we establish the equality
  of two   complicated expressions that depend on $\gamma$ only via the powers
  $(w+\gamma)^n$. In particular we can take $\gamma=z_2$, and the
  equalities will still hold, although the proof of Lemma
  \ref{lem:compatwithXS} breaks down in that case, as
  $S^{(\tau)}_{z_1,0}$ need not be defined.\label{rem:gammaz2valid}
\end{remark}

\begin{prop}\label{prop:ExpansionsX}
  Let $V$ be an $H_D$-quantum vertex algebra, and $A=a\ten b\ten c\in
  V^{\ten 3}$. Then we have the following expansions:
\begin{align}
    i_{z_1;z_2}X_{z_1,z_2,0}(A)&=Y(a,z_1)Y(b,z_2)c, \label{eq:propY1xY}\\
    i_{z_2;z_1}X_{z_1,z_2,0}(A)&=Y_{z_2}(1\otimes Y_{z_1})
    i_{z_2;z_1}S^{(\tau), 12}_{z_2,z_1}
    (b\otimes a\otimes c),\label{eq:propY1xYtau} \\
    i_{z_2;z_3}X_{z_2+z_3,z_2,0}(A)&=Y_{z_2}(Y_{z_3}\otimes 1)
    i_{z_2;z_3}S^{(z_2),12}_{z_3,0} (a\otimes b\otimes
    c).\label{eq:propYYx1 }
  \end{align}
\end{prop}

\begin{proof}
  \eqref{eq:propY1xY} is \eqref{eq:FullExpansX} for $n=3$ and $z_3=0$.
  By Corollary \ref{cor:permsymmetryX} (for $n=3$ and $f=\w_1$) we
  have
\[
X_{z_1,z_2,0}(A)=X_{z_2,z_1,0}\left (S^{(\tau)}_{z_2,z_1}(b\ten a)\ten c\right).
\]
Expanding this equation by applying $i_{z_2;z_1}$ and using
\eqref{eq:propY1xY} and definition \eqref{eq:DefY} gives \eqref{eq:propY1xYtau}.

 For the last part, let $\mathbf{f}=(132)=\w_2\w_1$, so that
 $\sigma_{\mathbf{f}}(a\ten b\ten c)=c\ten a\ten b$. Then
\begin{align*}
  S(\mathbf{f})&=S^{\mathbf{f}}_{z_1,z_2,z_3}\sigma_{\mathbf{f}}= 
  S^{\w_1}_{z_1,z_2,z_3}\mathbf{\tau}_1S^{\w_2}_{\w_1(z_1,z_2,z_3)}\tau_1=\\
  &=S^{(\tau)12}_{z_1,z_2}S^{(\tau)13}_{z_1,z_3}\sigma_{\mathbf{f}}.
\end{align*}
Therefore
\begin{equation}
  \label{eq:Sff}
  S^{\mathbf{f}}_{\mathbf{f}(z_2+w,z_2,z_1)}=
  S^{\mathbf{f}}_{z+1,z_2+w,z_2}=S^{(\tau)12}_{z_1,z_2+w} S^{(\tau)13}_{z_1,z_2}.
\end{equation}
Let $E=X_{z_2,z_1}(X_{{z_3},0}i_{z_2;{z_3}}S^{(z_2)}_{{z_3},0}\ten 1)(A)$.
Then
\begin{small}
  \begin{align*}
    E&=X_{z_1,z_2}S^{(\tau)}_{z_1,z_2}\tau(X_{{z_3},0}i_{z_2;{z_3}}S^{(z_2)}_{{z_3},0}\ten
    1)(A)= &&\text{by Lem. \ref{lem:X2braiding}}\\
    &=X_{z_1,z_2}S^{(\tau)}_{z_1,z_2}(1\ten
    X_{{z_3},0}i_{z_2;{z_3}}S^{(z_2)}_{{z_3},0})(c\ten
    a\ten b)= &&\\
    &=X_{z_1,z_2} \left(1\ten
      X_{{z_3},0}i_{z_2;{z_3}}S^{(z_2)}_{{z_3},0}\right)i_{z_1-z_2,z_2;{z_3}}
    S^{(\tau)12}_{z_1,z_2+{z_3}}S^{(\tau)13}_{z_1,z_2}\sigma_{\mathbf{f}}(A)= &&\text{by Lemma}\\
    &\qquad\qquad\qquad\text{  \ref{lem:compatwithXS} and Remark \ref{rem:gammaz2valid}}\\
    &=i_{z_1-z_2,z_2;{z_3}}\left
      (X_{z_1,z_2+{z_3},z_2}S^{(\tau)12}_{z_1,z_2+{z_3}}S^{(\tau)13}_{z_1,z_2}\right)
    \sigma_{\mathbf{f}}(A)= &&\text{by Thm. \ref{thm:analcontn}}\\
    &=i_{z_1-z_2,z_2;{z_3}}\left
      (X_{z_1,z_2+{z_3},z_2}S^{\mathbf{f}}_{\mathbf{f}
        (z_2+{z_3},z_2,z_1)}\right)
    \sigma_{\mathbf{f}}(A)= &&\text{by \eqref{eq:Sff}}\\
    &=i_{z_1-z_2,z_2;{z_3}}X_{z_2+{z_3},z_2,z_1}(A &&\text{by Thm.
      \ref{thm:analcontn}}
  \end{align*}
\end{small}
Putting $z_1=0$ proves then \eqref{eq:propYYx1 }.
\end{proof}

\section{The Braided Jacobi Identity}
\label{sec:braid-Jac-identity}

In one approach to the usual vertex algebras the Jacobi identity for
vertex operators is the basic identity, see e.g., \cite{MR2023933}. In
this section we derive the braided analog in our context of
$H_D$-quantum vertex algebras.

Introduce some more notation. If $f(z_1,z_2)\in
\mathbb{C}[[z_1,z_2]][z_1\inv,z_2\inv, (z_1-z_2)\inv]$ define the
difference of expansions of $f$ as
\begin{equation}
\delta(f(z_1,z_2))=(i_{z_1;z_2}-i_{z_2;z_1})(f(z_1,z_2)).\label{eq:defdelta}
\end{equation}
For instance,
\begin{equation}
  \label{eq:DiracDeltaDistr}
  \delta\left(\frac1{z_1-z_2}\right)=\delta(z_1,z_2)=\sum_{n\in \mathbb Z}
  z_1^nz_2^{-n-1}.
\end{equation}
This is the usual \emph{Dirac Delta Distribution}.

Recall that in this paper we are always expanding all expressions in
positive powers of $t$. For instance, if we write $\frac1{z_1-tz_2}$
we mean $\sum_{n\ge 0} (tz_2)^n/z_1^{n+1}$. Thus we have, for
instance,
\begin{equation}
  \label{eq:deltaz-tw}
\delta\left(\frac{1}{z_1-tz_2}\right)=0.
\end{equation}

\begin{lem}\label{lem:DiracDeltaIdent}For all $f(z_1,z_2,z_3)\in
    \mathbb{C}[[z_1,z_2,z_3]][z_1\inv,z_2\inv,z_3\inv]$ we have
    \begin{multline*}
      i_{z_1;z_2}\left(\delta(z_1-z_2,z_3)f(z_1,z_2,z_1-z_2)\right)-
      i_{z_2;z_1}\left(\delta(z_1-z_2,z_3)f(z_1,z_2,z_1-z_2)\right)=\\
      i_{z_2;z_3}\left(\delta(z_1,z_2+z_3)f(z_2+z_3,z_2,z_3)\right).
    \end{multline*}
  \end{lem}
  \begin{proof}
    See for example Proposition 2.3.26 in \cite{MR2023933}.
  \end{proof}
  \begin{defn}
We will write $a(z)$ for the 1-variable vertex operator $Y(a,z)$, the
field associated to $a\in V$.
  \end{defn}
 
 \begin{thm}(\textbf{Braided Jacobi Identity})\label{thm:braid-jacobi-ident}
    Let $V$ be an $H_D$-quantum vertex algebra. For all $a,b,c\in V$
    we have the identity:
    \begin{multline*}
      i_{z_1;z_2}\delta(z_1-z_2,z_3)a(z_1)b(z_2)c-
      i_{z_2;z_1}\delta(z_1-z_2,z_3)Y_{z_2}(1\otimes Y_{z_1})
      S^{(\tau),12}_{z_2,z_1}(b\otimes a\otimes c)\\
      =i_{z_2;z_3}\delta(z_1,z_2+z_3)Y_{z_2}(Y_{z_3}\otimes
      1)S^{(z_2),12}_{z_3,0}(a\otimes b\otimes c)
    \end{multline*}
  \end{thm}

  \begin{proof}
    Let $V^*$ be the dual of $V$, fix $v^* \in V^* $ and let
    $\langle\, ,\rangle$ be the pairing $V^*\otimes V\to
    k=\mathbb{C}[[t]]$. Then for all $A=a\otimes b\otimes c\in
    V^{\otimes 3}$ we have
\[
\langle v^*,
X_{z_1,z_2,0}(A)\rangle=\sum_{p\ge0}\sum_{l,m,n\in\mathbb{Z}} \frac{
  g_{l,m,n,p}(z_1,z_2)}{(z_1-z_2)^l z_1^m z_2^n} t^p,
\]
for $g_{l,m,n,p}(z_1,z_2)\in \mathbb{C}[[z_1,z_2]]$. (The sum over $l,
m,n$ is finite, for each $p$.) Define then
\[
F(z_1,z_2,z_3)=\sum_{p\ge0}\sum_{l,m,n\in\mathbb{Z}} \frac{
  g_{l,m,n,p}(z_1,z_2)}{z_3^l z_1^m z_2^n} t^p\in
\mathbb{C}[[z_1,z_2,z_3]][z_1\inv,z_2\inv,z_3\inv][[t]].
\]
Then we have by Corollary \ref{prop:ExpansionsX}
\begin{align*}
  i_{z_1;z_2}F(z_1,z_2,z_1-z_2)&=\langle v^*, a(z_1)b(z_2)c\rangle,\\
  i_{z_2;z_1}F(z_1,z_2,z_1-z_2)&=\langle v^*, Y_{z_2}(1\otimes
  Y_{z_1})i_{z_2;z_1}
  S^{(\tau),12}_{z_2,z_1}(b\otimes a\otimes c)\rangle,\\
  i_{z_2;z_3}F(z_2+z_3,z_2,z_3)&=\langle v^*, Y_{z_2}(Y_{z_3}\otimes
  1)i_{z_2;z_3}S^{(z_2),12}_{z_3,0}(a\otimes b\otimes c)\rangle.
\end{align*}
Then  we get from Lemma \ref{lem:DiracDeltaIdent} that
\begin{multline*}
  \langle v^*, i_{z_1;z_2}\delta(z_1-z_2,z_3)a(z_1)b(z_2)c\rangle-\\
 - \langle v^* ,i_{z_2;z_1}\delta(z_1-z_2,z_3)Y_{z_2}(1\otimes Y_{z_1})
  S^{(\tau),12}_{z_2,z_1}(b\otimes a\otimes c)\rangle\\
  =\langle v^*, i_{z_2;z_3}\delta(z_1,z_2+z_3)Y_{z_2}(Y_{z_3}\otimes
  1)i_{z_2,z_3}S^{(z_2),12}_{z_3,0}(a\otimes b\otimes c)\rangle.
   \end{multline*}
Since this is true for all $v^*\in V^*$ the Theorem follows.
  \end{proof}

  \begin{remark}
    Suppose $V$ is an $H_D$-quantum vertex algebra where
    $S^{(\tau)}_{z_1,z_2}$ and $S^{(\gamma)}_{z_1,z_2}$ both are the
    identity map on $V\otimes V$. Then the fields $a(z)=Y(a,z)$ satisfy
    the usual Jacobi identity:
\begin{multline*}
      i_{z_1;z_2}\delta(z_1-z_2,z_3)a(z_1)b(z_2)-
      i_{z_2;z_1}\delta(z_1-z_2,z_3)b(z_2)a(z_1)\\
      =i_{z_2;z_3}\delta(z_1,z_2+z_3)Y(Y(a,z_3)b,z_2),
    \end{multline*}
and it follows that $V$ is an ordinary vertex algebra, cf., \cite{MR2023933}. 

\end{remark}

 \section{Braided Borcherds Identity}
 \label{sec:BrBoRSlov}

 The original definition by Borcherds of vertex algebras was given
 in \cite{MR843307}. He took as starting point what later was called
 the Borcherds identity, instead of the Jacobi identity,
 cf., \cite{MR1651389}. In this section we derive a braided version of the
 Borcherds Identity.

 The following lemma is easy to check and well known (at least for
 $t=0$, see e.g., ??).
 \begin{lem}\label{lem:Tripleexpansionidentity}
   Let $W$ be a free $k$-module and $f(z,w)\in
   W[[z_1,z_2]][z_1^{-1},z_2^{-1},(z_1 - z_2)^{-1}][[t]]$. Then
 \[
 \Res_{z_1}\Big( \delta(f(z_1,z_2))\Big)=\Res_{z_3}\Big(i_{z_2;z_3}f(z_2+z_3,z_2)\Big).
 \]
 \end{lem}
 \begin{thm}(\textbf{Braided Borcherds Identity})\label{thm:BraidedBorcherds}
   Let $V$ be an $H_D$-quantum vertex algebra. Let $F\in\mathbb
   C[[z,w]][z^{-1},w^{-1},(z-w)^{-1}][[t]]$ and $a, b, c\in V$. Then
   we have the following identity:
 \begin{multline*}
   \Res_{z_1}\Big( Y(a,z_1)Y(b,z_2)c\,i_{z_1;z_2}F(z_1,z_2)-\\
   Y_{z_2}(1\otimes
   Y_{z_1})i_{z_2;z_1}S^{(\tau),12}_{z_2,z_1}(b\otimes a\otimes c)F(z_1,z_2)\Big)=\\
   = \Res_{z_3}\left( Y_{z_2}(Y_{z_3}\otimes 1)
     \,i_{z_2;z_3}\left(S^{(z_2), 12}_{z_3,0}(a\otimes b\otimes c)
       F(z_2+z_3,z_2)\right)\right).
 \end{multline*}
   \qed
 \end{thm}
 \begin{proof}
   Take in Lemma \ref{lem:Tripleexpansionidentity}
   $f(z_1,z_2)=X_{z_1,z_2,0}(a\otimes b\otimes c)F(z_1,z_2)$ and use
   Corollary \ref{prop:ExpansionsX} to relate expansions of $f$ to
   products and iterates of one-variable vertex operators.
 \end{proof}

\section{The $S$-Commutator,  Locality, and $(n)$-Products of Fields}
\label{sec:Scomloc}

\begin{defn}\label{defn:comS} Let $V$ be an $H_D$-quantum vertex
  algebra, and let $a,b,c\in V$. The \emph{$S$-commutator} of the
  fields associated to $a,b$ is
\[
\comS{a(z_1),b(z_2)}c=\delta\left(X_{z_1,z_2,0}(a\otimes b\otimes
  c)\right).
\]
\end{defn}

Here $\delta$ is the difference of expansions, see \eqref{eq:defdelta}.
We can write the $S$-commutator using Corollary \ref{prop:ExpansionsX}
explicitly as
\[
\comS{a(z_1),b(z_2)}c=
a(z_1)b(z_2)-Y_{z_2}(1\otimes Y_{z_1})i_{z_2;z_1} S_{z_2,z_1}^{(\tau)
 ,12}(b\otimes a\otimes c).
\]
Now the image of $\delta$ is a powerseries in $t$ with coefficients
(finite) sums of derivatives of the Dirac distribution
\eqref{eq:DiracDeltaDistr} with coefficients $V$-valued distributions
in $z_2$. So we can write the commutator as
\begin{equation}
  \label{eq:expscom}
\begin{aligned}
    \comS{a(z_1),b(z_2)}  &=\sum_{k>0} t^k
    \left(\sum_n\gamma_{n,k}(z_2)\partial_{z_2}^{(n)}\delta(z_1,z_2) \right)=\\
&=\sum_{n\ge0}
\gamma_n(z_2;t)\partial_{z_2}^{(n)}\delta(z_1,z_2).
\end{aligned}
  \end{equation}
This implies that for all $k\ge0$ there is an $N>0$ such that
\begin{equation}
  (z_1-z_2)^N \comS{a(z_1),b(z_2)}\equiv 0 \mod {t^k},
\label{eq:Slocality}
\end{equation}
and we see that the $S$-commutator of $a, b\in V$ is a local
distribution$\mod t^k$, see \cite{MR1651389}. (The $S$-commutator is
of course not necessarily itself local.)

\begin{defn}\label{defn:nproductfields}
  For all $n\in\mathbb{Z}$ the $(n)$-product of fields associated to
  $a,b\in V$ is
\[
a(z_2)_{(n)}b(z_2)c=\Res_{z_1}\big(\delta\Big(X_{z_1,z_2,0}(a\otimes
b\otimes c)(z_1-z_2)^n\Big)\big).
\]
\end{defn}

This definition allows us to write the $S$-commutator in terms of the
$(n)$-product of fields, for $n\ge0$.

\begin{thm}\label{thm:comsnproduct} Let $V$ be an $H_D$-quantum vertex algebra. For all $a,b\in V$
\[
\comS{a(z_1),b(z_2)}=\sum_{n\ge0}a(z_2)_{(n)}b(z_2)\partial_{z_2}^{(n)}\delta(z_1,z_2).
\]
\end{thm}
\begin{proof}
  By the usual calculus of local distributions, see e.g.,
  \cite{MR1651389}, it follows from \eqref{eq:expscom} that
  \begin{align*}
    \gamma_n(z_1;t)&=\Res_{z_2}\left(\comS{a(z_1),b(z_2)}(z_1-z_2)^n\right)=\\
  &=\Res_{z_2}\left( \delta\Big(X_{z_1,z_2,0}(a\otimes b\otimes -)(z_1-z_2)^n\Big)\right),
  \end{align*}
by Definition \ref{defn:comS}. Then the Lemma follows from Definition \ref{defn:nproductfields}.
\end{proof}

\section{$(n)$-Products of States}
\label{sec:n-productsstates}

We will call an element of $V$ also a \emph{state}. 
We define the $(n)$-product of states (as opposed to that of fields)
in $V$ in the usual way:
\[
a_{(n)}b=\Res_{z}\left (Y(a,z)bz^n \right),
\]
so that 
\begin{equation}
Y(a,z)=\sum_{n\in\mathbb Z}a_{(n)}z^{-n-1}.\label{eq:expansionY}
\end{equation}
We also have
\begin{equation}
  \label{eq:filedproperty}
  a_{(n)}b=0,\quad n\gg 0.
\end{equation}
In contrast to the usual vertex algebras the state-field
correspondence $a\mapsto a(z)$ is not quite a homomorphism of the
corresponding $(n)$-products: in general
\[
 a_{(n)}b(z)\ne a(z)_{(n)}b(z).
\]
Indeed, introduce the generating series ${\mathcal Y}_{{\mathcal F}}$ of the
$(n)$-products of fields by
\[
{\mathcal Y}_{{\mathcal F}}(a(z),w)=\sum_{n\in\mathbb{Z}} a(z)_{(n)}
w^{-n-1}.
\]
Then, if the state-field correspondence were a homomorphism we would
have
\begin{equation}
  \label{eq:StateFieldCorrHom}
\mathcal Y_{{\mathcal F}}(a(z),w)b(z)=Y(Y(a,w)b,z).  
\end{equation}
But this in general not true: the translation map
$S^{(\gamma)}_{z_1,z_2}$ is the obstruction to
(\ref{eq:StateFieldCorrHom}) being true. More precisely we have the
following theorem.

\begin{thm}\label{thm:StfieldCorr}
  \[
\mathcal Y_{{\mathcal F}}(a(z),w)b(z)c=Y_z(Y_w\ten 1)
i_{z;w}S^{(z)}_{w,0}(a\ten b)\ten c.  \]
\end{thm}

\begin{proof}
  By definition of the $(n)$-product of fields, Lemma
  \ref{lem:Tripleexpansionidentity} and Proposition
  \ref{prop:ExpansionsX} we have
\begin{align*}
  \mathcal Y_{{\mathcal
      F}}(a(z),w)b(z)c&=\Res_{z_1}\left(\delta(X_{z_1,z_2,0}(a\ten
    b\ten c)(z_1-z)^n)w^{-n-1}\right)=\\
  &=\Res_{z_3}\left(i_{z;z_3}(X_{z+z_3,z_2,0}(a\ten
    b\ten c)\delta(z_3,w)\right)=\\
    &=i_{z;w} X_{z+w,z,0}(a\ten b\ten c)\\
    &=Y_z(Y_w\ten 1)i_{z;w}S^{(z)}_{w,0}(a\ten b)\ten c.
\end{align*}
\end{proof}

Suppose that the translation map $S_{z_2,0}^{(z_3)}$ is such that
there exists $N\in\mathbb{Z}$ such that for all $a,b\in V$
\begin{equation}
i_{z_2;z_3} S^{(z_2)}_{z_3,0}(a\otimes b)=\sum_{k\ge -N}\left(\sum_{i}
a_{i,k}\otimes b_{i,k}\right) s_k(z_2) z_3^k,\quad s_k(z_2)\in \mathbb{C}((z_2)),
\label{eq:expansionTrans}
\end{equation}
where for fixed $k$ the summation over $i$ is finite.

Note that in a general $H_D$-quantum vertex algebra such expansion
need not exist. In the main example (see Section \ref{sec:mainexample})
this condition \emph{is}  satisfied, however.
\begin{cor}\label{cor:STFieldhom} Assume that
  \eqref{eq:expansionTrans} holds in $V$. Then for all $a,b\in V$ and
  $n\in \mathbb Z$
\[
a(z)_{(n)}b(z)
=\sum_{k\ge -N}\left(\sum_{i}
Y((a_{i,k})_{(k+n)}b_{i,k},z)\right) s_k(z).
\]
\end{cor}
\begin{proof}
  This is the case $F=(z_1-z_2)^n$ of the braided Borcherds identity,
  Theorem \ref{thm:BraidedBorcherds}. Indeed, in this case the LHS is
  just the $(n)$-product of the fields $a(z_2)$ and $b(z_2)$ acting on
  $c$, see Definition \ref{defn:nproductfields} and Corollary
  \ref{prop:ExpansionsX}. On the other hand the RHS of the braided
  Borcherds identity is in this case
\begin{align*}
   \Res_{z_3}&\left( 
 Y_{z_2} (Y_{z_3}\otimes 1)
 \sum_{k\ge -N}\left(\sum_{i} a_{i,k}\otimes b_{i,k}\otimes c\right) 
 s_k(z_2)z_3^{k+n}
\right)=
\\  
&= \Res_{z_3}\left( Y(\sum_{k\ge -N}\left(\sum_{i}
 (a_{i,k})_{(m)} b_{i,k}\right)z_3^{-m-1}, z_2) c
 s_k(z_2)z_3^{k+n}\right)=
\\
&=\sum_{k\ge -N}\left(\sum_{i}
 Y\left((a_{i,k})_{(k+n)} b_{i,k},z_2\right) c\right) s_k(z_2).
\end{align*}
The proof is concluded by the substitution $z_2\mapsto z$.
\end{proof}

\section{Normal Ordered Products and Operator Product Expansion}
\label{sec:NormOrdProd}

We have used the $(n)$-product (of fields) for $n\ge 0$ to calculate the
$S$-commutator, see Theorem \ref{thm:comsnproduct}. The $(n)$-products
for $n\le -1$ are also of course important.

\begin{defn}
  The normal ordered product of fields $a(z_1)$ and  $b(z_2)$ is given
  by
\[
\nopS{a(z_1)b(z_2)}c=\Res_z\left(\delta \Big(X_{z,z_2,0}(a\otimes
    b\otimes c)\frac1{z-z_1}\Big)\right).
\]
\end{defn}

We introduce projections on singular and holomorphic parts of a
formal distribution as usual by
\begin{align*}
  \Sing_{z_1}(f(z_1,z_2, \dots))&=-\Res_z\left( f(z,
    z_2,\dots)i_{z_1;z} \frac1{z-z_1}\right),\\
  \Hol_{z_1} (f(z_1,z_2, \dots))&= \Res_z\left( f(z,
    z_2,\dots)i_{z;z_1}\frac1{z-z_1}\right).
\end{align*}
In particular, if $f$ does not depend on $z_2, \dots$ we write
\[
  f_{\Sing}(z_1)=\Sing_{z_1}(f(z_1)),\quad
  f_{\Hol}(z_1)=\Hol_{z_1}(a(z_1)).
\]
Then we can rewrite the definition of the normal ordered product as
\[
\nopS{a(z_1)b(z_2)}=a_{\Hol}(z_1)b(z_2)+\Sing_{z_1}\Big(
Y_{z_2}(1\otimes Y_{z_1})i_{z_2;z_1}S^{(\tau)}_{z_2,z_1}(b\otimes
a)\Big).
\]
Comparing this with Definition \ref{defn:nproductfields} we see that
\[
a(z_2)_{(-1)}b(z_2)=\nopS{a(z_2)b(z_2)},
\]
and more generally
\[
a(z_2)_{(-n-1)}b(z_2)=\nopS{\partial_{z_2}^{(n)}a(z_2)b(z_2)}.
\]
This gives the \emph{Operator Product Expansion} of fields $a(z_1)$,
$b(z_2)$:
\begin{align*}
  a(z_1)b(z_2)&=\nopS{a(z_1)b(z_2)}+\Sing_{z_1}\left(\comS{a(z_1),b(z_2)}\right)=\\
             &=\nopS{a(z_1)b(z_2)} +\sum_{n\ge0}a(z_2)_{(n)}b(z_2)
    i_{z_1;z_2}\left(\frac1{(z_1-z_2)^{n+1}}\right).
\end{align*}
Of coure, using Corollary \ref{cor:STFieldhom} we can express the operator
product expansion in terms of the $(n)$-product of states, but this
seems rather messy.

\section{Weak Associativity}
\label{sec:WeakAssoc}
Two basic ingredients in the usual theory of vertex algebras are
locality and associativity. For $H_D$-quantum vertex algebras the
analog of locality is $S$-locality, (\ref{eq:Slocality}). In this section we
derive the analog of associativity. It involves the translation map $S^{(\gamma)}_{z_1,z_2}$.
\begin{thm}[\textbf{Weak associativity}]
  Let $V$ be an $H_D$-quantum vertex algebra. For all $a,b, c\in V$
  and for all powers $t^k$ there is an $N\ge 0$ such that
  \begin{multline*}
(z_2+z_3)^N i_{z_3;z_2}a(z_2+z_3)b(z_2)c\equiv\\
\equiv (z_2+z_3)^N Y_{z_2}(Y_{z_3}\otimes
  1)i_{z_2;z_3}(S^{(z_2)}_{z_3,0}(a\otimes b)\otimes c) \mod t^k.
  \end{multline*}
\end{thm}

\begin{proof}
  Take $\Res_{z_1}$ in the braided Jacobi identity of Theorem
  \ref{thm:braid-jacobi-ident} to find
  \begin{align*}
    i_{z_3;z_2}&a(z_2+z_3)b(z_2)c-
    Y_{z_2} (Y_{z_3}\otimes 1)
    (i_{z_2;z_3}S^{(z_2)}_{z_3,0}(a\otimes b)\otimes c)=\\
&=-\Res_{z_1}\left(i_{z_2;z_1}\delta(z_1-z_2,z_3)
  Y_{z_2}(1\otimes Y_{z_1}) (S^{(\tau)}_{z_2,z_1}(b\otimes a)\otimes c)\right)=\\
&=-\Res_{z_1}\left(\sum_{k=0}^\infty(-z_1)^k \partial_{z_2}^{(k)}\delta(-z_2,z_3)
  Y_{z_2}(1\otimes Y_{z_1})i_{z_2;z_1}(S^{(\tau)}_{z_2,z_1}(b\otimes a)\otimes c)\right).
  \end{align*}
Expanding the RHS observe that the coefficient of each power of $t$ is
after taking the residue a finite sum of $z_2$ derivatives of
$\delta(-z_2,z_3)$, hence vanishes if multiplied by a suitable power
of $z_2+z_3$.
\end{proof}
\begin{remark}
  For ordinary vertex algebras the power of $N$ in weak associativity
  depends only on $a$ and $c$, not on $b$. The above proof in the case
  of $H_D$-quantum vertex algebras does not allow us to conclude the
  same, because of the appearance of the braiding
  $S^{(\tau)}_{z_2,z_1}(b\ten a)$.
\end{remark}

\section{The $H_D$-Bialgebra $V$}
\label{sec:BialgV}

In the rest of the paper we will construct a class of examples of
$H_D$-quantum vertex algebras, using bicharacters on the underlying
space $V$. To define bicharacters we need to assume that $V$ has extra
structure: we will assume that $V$ is a commutative and cocommutative
$k$-bialgebra, or even a Hopf algebra. The coproduct and counit of $V$
will be denoted by $\Delta$ and $\epsilon$. We assume also that $V$ has
a compatible $H_D$-action. This means that
\begin{itemize}
\item $D(ab)=(Da)b+aDb$, $a, b\in V$.
\item $\Delta(Da)=\Delta_{H_D}(D)\Delta(a)$, $a\in V$.
\item $\epsilon(Da)=\epsilon_{H_D}(D)\epsilon(a)=0$.
\end{itemize}
We will call a $V$ as above an $H_D$-bialgebra. The identity element
$1=1_V$ will be the vacuum of $V$.

\section{Bicharacters}
\label{sec:BiChar}
Let $W_2$ be the algebra of power series in $t$, with coefficients
rational functions in $z_1,z_2$ with poles at $z_1=0$, $z_2=0$ or
$z_1=z_2$:
\begin{equation}
  \label{eq:DefW2}
W_2=\mathbb C[z_1^{\pm 1},z_2^{\pm 1}, (z_1-z_2)^{\pm 1}][[t]].
\end{equation}
We extract some results from \cite{MR1865087} on \emph{bicharacters}.
A $W_2$-valued bicharacter on an $H_D$-bialgebra $V$ is a linear map
\[
r_{z_1,z_2}\colon V^{\otimes 2}\to W_2,
\]
satisfying
\begin{itemize}
\item (\textbf{Vacuum}) $r_{z_1,z_2}(a\otimes 1)=r_{z_1,z_2}(1\otimes a)=\epsilon(a)$,
  $a\in V$.
\item (\textbf{Multiplication}) For all $a,b,c\in V$ we have
  $r_{z_1,z_2}(a\otimes bc)=\sum r_{z_1,z_2}(a^\prime\otimes
  b)r_{z_1,z_2}(a^{\prime\prime}\otimes c)$ and $r_{z_1,z_2}(ab\otimes
  c)=\sum r_{z_1,z_2}(a\otimes c^\prime)r_{z_1,z_2}(b\otimes
  c^{\prime\prime})$.
\end{itemize}
Here and below we use the notation $\Delta(a)=\sum a^\prime\otimes
a^{\prime\prime}$ for the coproduct of $a\in V$. Often we will also omit the
summation symbol, to unclutter the formulas.

In case the bicharacter additionally satisfies
\begin{itemize}
 \item (\textbf{$H_D\otimes H_D$-covariance}) $r_{z_1,z_2}(D^ka\otimes D^\ell
  b)=\partial_{z_1}^k\partial_{z_2}^\ell r_{z_1,z_2}(a\otimes b)$,
  $a,b\in V$,
\end{itemize}
we call the bicharacter \emph{$H_D\otimes H_D$-covariant}.

We can multiply bicharacters: 
\begin{equation}
(r\ast s)_{z_1,z_2}(a\otimes b)=r_{z_1,z_2}(a^\prime\otimes b^\prime)
s_{z_1,z_2}(a^{\prime\prime}\otimes b^{\prime\prime
}).\label{eq:multbichar}
\end{equation}
The unit bicharacter is 
\begin{equation}
\epsilon_{z_1,z_2}(a\otimes
b)=\epsilon(a)\epsilon(b).\label{eq:unitbichar}
\end{equation}
The collection of bicharacters on an $H_D$-bialgebra forms then a
commutative mo\-no\-id. 

In case $V$ is an $H_D$-Hopf algebra, i.e., comes with an antipode
compatible with the $H_D$-action, all bicharacters are invertible,
with inverse given by
\[
r\inv_{z_1,z_2}(a\otimes b)=r_{z_1,z_2}(S(a)\otimes b).
\]
In this case the set of bicharacters forms an Abelian group.

The transpose of a bicharacter is defined by
\[
r^\tau_{z_1,z_2}(a\otimes b)=r_{z_2,z_1}(b\otimes a).
\]
The transpose is an involution of the monoid of bicharacters:
\[
(r\ast s)_{z_1,z_2}^\tau =(r^\tau\ast s^\tau)_{z_1,z_2}.
\]
If $r$ is an invertible bicharacter with inverse $r\inv$ we relate the
transpose $r^\tau$ to $r$ by
\begin{equation}
  \label{eq:braidingforr}
r^\tau_{z_1,z_2}=r_{z_1,z_2}\ast R_{z_1,z_2},
\end{equation}
where
\begin{equation}
  \label{eq:defR}
R_{z_1,z_2}=r\inv_{z_1,z_2}\ast r^\tau_{z_1,z_2}.
\end{equation}
We will call $R_{z_1,z_2}$ the \emph{braiding bicharacter} associated
to $r_{z_1,z_2}$. It is the obstruction to $r$ being \emph{symmetric}:
$r=r^\tau$. It will control the braiding in the quantum vertex algebra
we are going to construct from $r_{z_1,z_2}$ in Section
\ref{sec:QuantVertBichar} below.  The braiding bicharacter
$R_{z_1,z_2}$ is \emph{unitary}:
\begin{equation}
  \label{eq:Runitary}
  R^\tau_{z_1,z_2}=R\inv_{z_1,z_2}.
\end{equation}
Define for a bicharacter $r_{z_1,z_2}$ a shift
\begin{equation}
  \label{eq:gammabichar}
  r^\gamma_{z_1,z_2}=r_{z_1+\gamma,z_2+\gamma}.
\end{equation}
The shift $r^\gamma_{z_1,z_2}$ is again a bicharacter. If
$r_{z_1,z_2}$ is $H_D\otimes H_D$-covariant we have the following
expansion:
\[
i_{z_1,z_2;\gamma} r^\gamma_{z_1,z_2}=r_{z_1,z_2}\circ \Delta(e^{\gamma D}).
\]
In case the bicharacter is invertible we relate the shift $r^\gamma$
to $r$ by
\begin{equation}
  \label{eq:DefRgamma}
r^\gamma_{z_1,z_2}=r_{z_1,z_2}\ast R^\gamma_{z_1,z_2}, \quad
R^\gamma_{z_1,z_2}=r\inv_{z_1,z_2}\ast r^\gamma_{z_1,z_2}.
\end{equation}
We call $R^\gamma_{z_1,z_2}$ the \emph{translation bicharacter}
associated to $r_{z_1,z_2}$. It is the obstruction to $r$ being shift
invariant (i.e., to $r$ being a function just of $z_1-z_2$).



\section{$H_D$-Quantum Vertex Algebras from Bicharacters}
\label{sec:QuantVertBichar}

Suppose now that $V$ is an $H_D$-bialgebra with invertible bicharacter
$r_{z_1,z_2}$. In general, a bicharacter on $V$ takes values in $W_2$,
see (\ref{eq:DefW2}). For the purpose of the construction of vertex
operators we need to make an extra assumption: that $r_{z_1,z_2}$ can
be evaluated at $z_2=0$. More precisely, we make the following
\begin{defn}\label{defn:VO-assumption} A bicharacter $r_{z_1,z_2}$
  satisfies the \emph{Vertex Operator Assumption} if it is a map
\begin{equation}
 \label{eq:VOassumption}  r_{z_1,z_2}\colon V^{\otimes 2}\to 
  \mathbb{C}[z_1^{\pm1},z_2,(z_1-z_2)\inv][[t]]. \tag{{VO assumption}}
\end{equation}
  \end{defn}
In the sequel we will use $\rho_{z_1,z_2}$ to denote an arbitrary
$W_2$-valued bicharacter, and we will write $r_{z_1,z_2}$ for a
bicharacter satisfying the \ref{eq:VOassumption}.

Following the general philosophy of Borcherds, \cite{MR1865087}, (but
not the technical details) we define in this section, given an
invertible bicharacter $r_{z_1,z_2}$ satisfying the
\ref{eq:VOassumption}, an $H_D$-quantum vertex algebra structure on
$V$. The final result is summarized in Theorem
\ref{thm:h_d-quantum-vertex-from-bichar} below,

We define for any bicharacter $\rho_{z_1,z_2}$ on $V$ a map
$S^{\rho_{z_1,z_2}}$ on $V\otimes V$ by
\begin{equation}
  \label{eq:defS(rho)}
  S^{\rho_{z_1,z_2}}(a\otimes b)=a^\prime\otimes b^\prime
  \rho_{z_1,z_2}(a^{\prime\prime}\otimes b^{\prime\prime}).
\end{equation}
In particular, to a bicharacter $r_{z_1,z_2}$ satisfying the
\ref{eq:VOassumption} with braiding bicharacter $R_{z_1,z_2}$, see
\eqref{eq:defR}, we associate the map
\begin{equation}
  \label{eq:defStau}
  S^{(\tau)}_{z_1,z_2}=S^{R_{z_1,z_2}}\colon V\otimes V\to   V\otimes
  V[z_1^{\pm1},z_2^{\pm1},(z_1-z_2)^{\pm1}][[t]], 
\end{equation}
and associated to the translation bicharacter \eqref{eq:DefRgamma} we
get a map
\begin{equation}
  \label{eq:defSgama}
  S^{(\gamma)}_{z_1,z_2}=S^{R^\gamma_{z_1,z_2}}\colon V\otimes V\to   V\otimes
  V[z_1^{\pm1},z_2,(z_1+\gamma)^{\pm1}, (z_2+\gamma), (z_1-z_2)^{\pm1}][[t]], 
\end{equation}

\begin{lem}\label{lem:Srhosigma}  
  \begin{enumerate}
  \item If $\epsilon$ is the unit bicharacter on $V$, then
  $S^{\epsilon}=1_{V^{\otimes2}}$.
\item If $\rho_{z_1,z_2},\sigma_{z_1,z_2}$ are bicharacters on $V$, then
  $S^{\rho_{z_1,z_2}\ast\sigma_{z_1,z_2}}=S^{\rho_{z_1,z_2}}\circ S^{\sigma_{z_1,z_2}}$.
\item If $\rho_{z_1,z_2}$ is a bicharacter, then $\tau\circ
  S^{\rho_{z_1,z_2}}\circ\tau=S^{\rho_{z_1,z_2}^\tau}$.
  \end{enumerate}
\end{lem}

Define then, for given invertible bicharacter $r_{z_1,z_2}$ satisfying
the \ref{eq:VOassumption}, singular multiplication maps
\[
X_{z_1,z_2}\colon V^{\otimes 2}\to
V\otimes[[z_1,z_2]][z_1\inv,(z_1-z_2)\inv][[t]]. 
\]
by
\begin{equation}
  \label{eq:DefSingmult}
  X_{z_1,z_2}=m_2\circ (e^{z_1D}\otimes e^{z_2D})\circ
  S^{r_{z_1,z_2}},
\end{equation}
where $m_2$ is the (nonsingular) multiplications of the (associative)
algebra $V$. More explicitly (dropping here and below the nonsingular
multiplication $m_2$ on $V$):
\[
X_{z_1,z_2}(a\otimes b)=e^{z_1D}a^\prime e^{z_2D}b^\prime
r_{z_1,z_2}(a^{\prime\prime}\otimes b^{\prime\prime}).
\]
\begin{lem}
  For any bicharacter $\rho_{z_1,z_2}$ on $V$ we have for $a\in V$
\[
S^{\rho_{z_1,z_2}}(a\otimes 1)=a\otimes 1,\quad S^{\rho_{z_1,z_2}}(1\otimes a)=a\otimes 1.
\]
\end{lem}
\begin{proof}
  Since $\rho_{z_1,z_2}$ is a bicharacter we have
\[
\rho_{z_1,z_2}(a\otimes 1)=\epsilon(a)=\rho_{z_1,z_2}(1\otimes a).
\]
In any bialgebra we have $a^\prime \epsilon(a^{\prime\prime})=a$, and
the Lemma follows from the definition of $S^{\rho_{z_1,z_2}}$, see \eqref{eq:defS(rho)}.
\end{proof}
\begin{cor}
  The vacuum axioms \eqref{eq:vacuumX} and \eqref{eq:vacuumS} hold for
  $X_{z_1,z_2}$ defined by \eqref{eq:DefSingmult} and for
  $S^{(\tau)}_{z_1,z_2}, S^{(\gamma)}_{z_1,z_2}$ defined by
  \eqref{eq:defStau} and \eqref{eq:defSgama}.
\end{cor}

\begin{lem}
  For any $H_D\otimes H_D$-covariant bicharacter $\rho_{z_1,z_2}$ we
  have
\[
[S^{\rho_{z_1,z_2}},1\otimes D]=\partial_{z_2}S^{\rho_{z_1,z_2}},\quad
[S^{\rho_{z_1,z_2}},D\otimes 1]=\partial_{z_1}S^{\rho_{z_1,z_2}}.
\]
\end{lem}
\begin{proof}
  By assumption on $V$ we have  $\Delta(Db)=D b^\prime\otimes
  b^{\prime\prime}+ b^\prime\otimes D b^{\prime\prime}$. By assumption
  on the bicharacter we have $\rho_{z_1,z_2}(a\otimes
  Db)=\partial_{z_2}\rho_{z_1,z_2}(a\otimes b)$. Then, for $a,b\in V$
  \begin{align*}
    S^{\rho_{z_1,z_2}}(a\otimes Db)&=a^\prime\otimes D
    b^\prime\rho_{z_1,z_2}(a^{\prime\prime}\otimes b^{\prime\prime})+
    a^\prime\otimes b^\prime\rho_{z_1,z_2}(a^{\prime\prime}\otimes
    Db^{\prime\prime})=\\
    &=(1\otimes
    D)S^{\rho_{z_1,z_2}}+\partial_{z_2}S^{\rho_{z_1,z_2}}(a\otimes b),
  \end{align*}
proving the first part. The second part is similar.

\end{proof}

\begin{cor} The $H_D$-covariance axiom \eqref{eq:H_DcovX} holds for
  $X_{z_1,z_2}$ defined by \eqref{eq:DefSingmult} and  the $H_D$-covariance axiom \eqref{eq:H_DcovS} holds for
  $S^{(\tau)}_{z_1,z_2}, S^{(\gamma)}_{z_1,z_2}$ defined by
  \eqref{eq:defStau} and \eqref{eq:defSgama}.
\end{cor}

\begin{lem}
  The $H_D$-covariance axiom (\ref{eq:H_DcovMult}) holds for
  $X_{z_1,z_2}$ defined by (\ref{eq:DefSingmult}).
\end{lem}

\begin{proof}
We have
\begin{align*}
  X_{z_1 + \gamma,z_2+\gamma} (a\ten b)  &=e^{(z_1+\gamma)D}a^\prime
  e^{(z_2+\gamma)D}b^\prime
  r_{z_1+\gamma,z_2+\gamma}(a^{\prime\prime}\ten b^{\prime\prime})=\\
&=e^{\gamma D}(e^{z_1D}a^\prime  e^{(z_2)D}b^\prime)
  r_{z_1,z_2}({a^{\prime\prime}}^\prime\ten
  {b^{\prime\prime}}^\prime)R^\gamma_{z_1,z_2}({a^{\prime\prime}}^{\prime\prime}
\ten {b^{\prime\prime}}^{\prime\prime})=\\
&e^{\gamma D} X_{z_1,z_2}\circ S^{(\gamma)}_{z_1,z_2}(a\ten b).
\end{align*}
\end{proof}

\begin{lem}
  For any bicharacter $\rho_{z_1,z_2}$ the map $S^{\rho_{z_1,z_2}}$
  satisfies the Yang-Baxter equation (\ref{eq:YBaxiom}).
\end{lem}

\begin{proof}
  This follows form the combined cocommutativity and coassociativity
  identity
\[
\tau^{23}(\Delta\otimes 1)\Delta=(\Delta\otimes 1)\Delta.
\]
  \end{proof}

  \begin{cor}
The maps
  $S^{(\tau)}_{z_1,z_2}, S^{(\gamma)}_{z_1,z_2}$ defined by
  \eqref{eq:defStau} and \eqref{eq:defSgama} satisfy the Yang-Baxter
  axiom \eqref{eq:YBaxiom}.    
  \end{cor}

  \begin{lem}
    For any bicharacter $\rho_{z_1,z_2}$ the map $S^{\rho_{z_1,z_2}}$
    is compatible with the singular multiplication:
    \begin{align*}
      S^{\rho_{z_1,z_2}}(X_{w_1,w_2}\otimes 1)&=(X_{w_1,w_2}\otimes 1)
      i_{z_1,z_1-z_2;w_1,w_2}S^{\rho_{z_1+w_1,z_2},23}S^{\rho_{z_1+w_2,z_2},13},\\
      S^{\rho_{z_1,z_2}}(1\otimes X_{w_1,w_2})&=(1\otimes X_{w_1,w_2})
      i_{z_1-z_2,z_2;w_1,w_2}S^{\rho_{z_1,z_2+w_1},12}S^{\rho_{z_1,z_2+w_2},13}.
    \end{align*}\end{lem}

\begin{proof}
  For $a,b,c\in V$ we have
  \begin{align*}
    S^{\rho_{z_1,z_2}}(&X_{w_1,w_2}\otimes 1)(a\otimes b\otimes c)=
    S^{\rho_{z_1,z_2}}(e^{w_1D}a^\prime e^{w_2D} b^\prime\otimes
    c)r_{w_1,w_2}(a^{\prime\prime}\otimes
    b^{\prime\prime})=\\
    &=(e^{w_1D}a^\prime e^{w_2D} b^\prime)^\prime
    c^\prime\rho_{z_1,z_2}((e^{w_1D}a^\prime e^{w_2D}
    b^\prime)^{\prime\prime}\otimes c^{\prime
      \prime})r_{w_1,w_2}(a^{\prime\prime}\otimes
    b^{\prime\prime})=\\
     &= e^{w_1D} {a^\prime}^\prime e^{w_2D} {b^\prime}^\prime
     c^\prime i_{z_1,z_1-z_2;w_1,w_2}
     \rho_{z_1+w_1,z_2} ({a^\prime}^{\prime\prime}\otimes
     {c^{\prime\prime}}^\prime) 
     \rho_{z_1+w_2,z_2} ({b^\prime}^{\prime\prime}\otimes {c^{\prime
       \prime}}^{\prime \prime})\times\\
     & \qquad\qquad\qquad\qquad\qquad\qquad\qquad\qquad \times
     r_{w_1,w_2}(a^{\prime\prime}\otimes b^{\prime\prime})=\\
     &=i_{z_1,z_1-z_2;w_1,w_2}(X_{w_1,w_2}\otimes
      1)S^{\rho_{z_1+w_1,z_2},23}S^{\rho_{z_1+w_2,z_2},13}(a\otimes
      b\otimes c)
\end{align*}
The proof of the other part is similar.
\end{proof}

\begin{cor}
  $S^{(\tau)}_{z_1,z_2}, S^{(\gamma)}_{z_1,z_2}$ defined by
  \eqref{eq:defStau} and \eqref{eq:defSgama} satisfy the compatibility
  with multiplication axioms (\ref{eq:CompatXx1}) and (\ref{eq:Compat1xX}).
\end{cor}

\begin{cor}
$S^{(\tau)}_{z_1,z_2}$ defined by
  \eqref{eq:defStau}  satisfies the unitarity axiom (\ref{eq:GrpSStau}).
  \end{cor}
  \begin{proof}
For unitarity, recall that $S^{(\tau)}_{z_1,z_2}=S^{R_{z_1,z_2}}$, so
that by Lemma \ref{lem:Srhosigma} we have
  $\tau\circ S^{(\tau)}_{z_1,z_2}\circ \tau=S^{R_{z_1,z_2}^\tau}$ so that by Lemma
  \ref{lem:Srhosigma} again and (\ref{eq:Runitary}) we find
  \[
  S^{(\tau)}_{z_1,z_2}\circ \tau \circ S^{(\tau)}_{z_1,z_2}\circ
  \tau=S^{R_{z_1,z_2}}\circ S^{R^\tau_{z_1,z_2}}=
  S^{\epsilon}=1_{V\otimes V}.
  \]    
  \end{proof}

  \begin{cor}
    $S^{(\gamma)}_{z_1,z_2}$ defined by \eqref{eq:defSgama} satisfies
    the group axioms (\ref{eq:SatZero}) and (\ref{eq:GrpSgam1gam2}).
  \end{cor}

  \begin{proof}
    Since $R_{z_1,z_2}^{(\gamma=0)}=\epsilon$, the unit bicharacter,
    axiom (\ref{eq:SatZero}) follows.

    Now
    \begin{align*}
      R^{\gamma_1+\gamma_2}_{z_1,z_2}&=r\inv_{z_1,z_2}\ast
      r^{\gamma_1+\gamma_2}_{z_1,z_2}=\\
      &= r\inv_{z_1,z_2}\ast
      r^{\gamma_1}_{z_1,z_2}\ast r^{\gamma_1, -1}_{z_1,z_2}\ast
      r^{\gamma_1+\gamma_2}_{z_1,z_2}=\\
      &=R^{\gamma_1}_{z_1,z_2} \ast R^{\gamma_2}_{z_1+\gamma_1,z_2+\gamma_1}, 
    \end{align*}
    so that axiom (\ref{eq:GrpSgam1gam2}) follows from Lemma
    \ref{lem:Srhosigma} and definition (\ref{eq:defSgama}).
  \end{proof}

  \begin{lem}
    $X_{z_1,z_2}$ and $S^{(z_2)}_{w_1,w_2}$
    defined by (\ref{eq:DefSingmult}) and (\ref{eq:defSgama}) satisfy
    the locality Axiom \eqref{eq:localityAx}.
  \end{lem}
  \begin{proof}
Define
\[
E=e^{z_1D}a^\prime e^{z_2 D}b^\prime c^\prime
r_{z_1,z_2}({a^{\prime\prime}}^\prime\ten{b^{\prime\prime}}^\prime)
r_{z_1,0}({a^{\prime\prime}}^{\prime\prime}\ten{c^{\prime\prime}}^\prime)
r_{z_2,0}({b^{\prime\prime}}^{\prime\prime}\ten{c^{\prime\prime}}^{\prime\prime}).
\] 
Then
\begin{align*}
  X_{z_1,0} &(1\ten X_{z_2,0})(A)
  =e^{z_1D}a^\prime
  \left(e^{z_2D}b^\prime c^\prime\right)^\prime 
  r_{z_1,0}(a^{\prime\prime}\ten  \left(e^{z_2D}b^\prime
    c^\prime\right)^{\prime\prime })
  r_{z_2,0}(b^{\prime\prime}\ten c^{\prime\prime})=\\
&=e^{z_1D}a^\prime e^{z_2D}{b^\prime}^\prime {c^\prime}^\prime 
  r_{z_1,0}({a^{\prime\prime}}^\prime\ten  e^{z_2D}{b^\prime}^{\prime\prime})
  r_{z_1,0}({a^{\prime\prime}}^{\prime\prime}\ten  {c^\prime}^{\prime\prime})
  r_{z_2,0}(b^{\prime\prime}\ten c^{\prime\prime})=\\
&=e^{z_1D}a^\prime e^{z_2D}{b^\prime}{c^\prime}^\prime 
  i_{z_1;z_2}r_{z_1,z_2} ({a^{\prime\prime}}^\prime\ten  {b^{\prime\prime}}^{\prime})
  r_{z_1,0} ({a^{\prime\prime}}^{\prime\prime}\ten  {c^{\prime\prime}}^{\prime})
  r_{z_2,0}({b^{\prime\prime}}^{\prime\prime}\ten
  {c^{\prime\prime}}^{\prime\prime})=\\
&=i_{z_1;z_2}E.
\end{align*}
On the other hand
\begin{align*}
  X_{z_2,0} &(1\ten X_{z_1,0})i_{z_2;z_1}S^{(\tau)12}_{z_2,z_1}(b\ten
  a\ten c)=\\
  &=X_{z_2,0} (1\ten X_{z_1,0})(b^\prime\ten a^\prime \ten
  c)i_{z_2;z_1}
  R^{(\tau)}_{z_2,z_1}(b^{\prime\prime}\ten a^{\prime\prime})=\\
  &=X_{z_2,0} (b^\prime \ten e^{z_1D} {a^\prime}^\prime c^\prime)
  r_{z_1,0}( {a^\prime}^{\prime\prime} \ten c^{\prime\prime})
  i_{z_2;z_1} R^{(\tau)}_{z_2,z_1}(b^{\prime\prime}\ten a^{\prime\prime})=\\
  &=e^{z_2D} {b^\prime}^\prime (e^{z_1D} {a^\prime}^\prime
  c^\prime)^\prime r_{z_2,0}( {b^\prime}^{\prime\prime}\ten (e^{z_1D}
  {a^\prime}^\prime
  c^\prime)^{\prime\prime} )\times\\
  &\qquad\quad \times r_{z_1,0}( {a^\prime}^{\prime\prime} \ten
  c^{\prime\prime}) i_{z_2;z_1}
  R^{(\tau)}_{z_2,z_1}(b^{\prime\prime}\ten
  a^{\prime\prime})=\\
  &=e^{z_2D}b^\prime e^{z_1D} a^\prime c^\prime r_{z_2,0}
  ({{b^{\prime\prime}}^{\prime}}^\prime \ten e^{z_1D}
  {{a^{\prime\prime}}^\prime}^\prime)
  r_{z_2,0}({b^{\prime\prime}}^{\prime\prime}\ten  {c^{\prime\prime}}^{\prime})\times\\
  &\qquad\quad \times r_{z_1,0}( {a^{\prime\prime}}^{\prime\prime}
  \ten {c^{\prime\prime}}^{\prime\prime}) i_{z_2;z_1}
  R^{(\tau)}_{z_2,z_1}({{b^{\prime\prime}}^\prime}^{\prime\prime}\ten
  {{a^{\prime\prime}}^\prime}^{\prime\prime})=\\
  &=e^{z_2D}b^\prime e^{z_1D} a^\prime c^\prime
  i_{z_2,z_1}\left(r_{z_2,z_1}({{b^{\prime\prime}}^{\prime}}^\prime
    \ten {{a^{\prime\prime}}^\prime}^\prime)
    R^{(\tau)}_{z_2,z_1}({{b^{\prime\prime}}^\prime}^{\prime\prime}\ten
    {{a^{\prime\prime}}^\prime}^{\prime\prime})\right)\times\\
  &\qquad\quad \times  r_{z_2,0}({b^{\prime\prime}}^{\prime\prime}\ten  {c^{\prime\prime}}^{\prime})
 r_{z_1,0}( {a^{\prime\prime}}^{\prime\prime}
  \ten {c^{\prime\prime}}^{\prime\prime}) =\\
  &=e^{z_2D}b^\prime e^{z_1D} a^\prime c^\prime
  i_{z_2,z_1}\left(r_{z_1,z_2}({{a^{\prime\prime}}^\prime}\ten
    {{b^{\prime\prime}}^{\prime}} )\right)
  r_{z_2,0}({b^{\prime\prime}}^{\prime\prime}\ten  {c^{\prime\prime}}^{\prime})\times\\
  &\qquad\quad \times r_{z_1,0}( {a^{\prime\prime}}^{\prime\prime}
  \ten {c^{\prime\prime}}^{\prime\prime}) =\\
  &=i_{z_2,z_1}E.
\end{align*}
Since
\[
(z_1-z_2)^N i_{z_1;z_2}E=(z_1-z_2)^N i_{z_2;z_1}E
\]
the locality Axiom \eqref{eq:localityAx} follows.
  \end{proof}

The results in this section are summarized in the following theorem.
  \begin{thm}\label{thm:h_d-quantum-vertex-from-bichar}
    Let $V$ be an $H_D$-bialgebra with invertible bicharacter
    $r_{z_1,z_2}$, satisfying the \ref{eq:VOassumption} of Definition
    \ref{defn:VO-assumption}. Then the singular multiplications
    $X_{z_1,z_2}$, $X_{z_1,z_2,z_3}$ and maps $S^{(\tau)}_{z_1,z_2}$,
    $S^{(\gamma)}_{z_1,z_2}$ defined by (\ref{eq:DefSingmult}),
    (\ref{eq:defStau}) and (\ref{eq:defSgama}) give $V$ the structure
    of an $H_D$-quantum vertex algebra as in Definition
    \ref{defn:h_d-quantum-vertex-alg} .
  \end{thm}

\section{Bicharacters and EK-quantum Vertex Operator Algebras}
\label{sec:BicharEKqva}

Let $V$ be an $H_D$-bialgebra with invertible bicharacter, so that we
have on $V$ by Theorem \ref{thm:h_d-quantum-vertex-from-bichar} an
$H_D$-quantum vertex algebra structure. In case the bicharacter satisfies
\begin{equation}
  \label{eq:Ginvr}
  (\partial_{z_1}+\partial_{z_2})r_{z_1,z_2}=0
\end{equation}
the bicharacter is really just a function of $z_1-z_2$: $r_{z_1,z_2}$ takes
values in $\mathbb C[(z_1-z_2)^{\pm 1}][[t]]$. In this case the
translation bicharacter $R^\gamma_{z_1,z_2}$ is the unit bicharacter on
$V$. 

In this situation we can evaluate the bicharacter $r_{z_1,z_2}$, the
vertex operator $X_{z_1,z_2}$ and the braiding $S_{z_1,z_2}$ both at
$z_1=0$ and at $z_2=0$.

We have in this case $r_{0,z}=r_{-z,0}$ so that
\[
X_{0,z}(a\otimes b)=e^{zD}\left(e^{-zD}a^\prime b^\prime\right)
r_{0,z}(a^{\prime\prime}\otimes b^{\prime\prime})=e^{zD}Y(a, -z)b.
\]
The braided commutativity Lemma \ref{lem:X2braiding} gives, by putting
$z_2=0$,
\[
Y(a^\prime,z)b^\prime R_{z,0}(a^{\prime\prime}\otimes
b^{\prime\prime})= e^{zD}Y(b,-z)a.
\]
We emphasize that in general $H_D$-quantum vertex algebras one does not
have a similar braided skew-symmetry, since the braiding
$S^{(\tau)}_{z_1,z_2}$ cannot be evaluated at $z_2=0$.

The $H_D$-covariance axiom \eqref{eq:H_DcovMult} reduces to
the familiar formula
\begin{equation}
  \label{eq:TransCovYclass}
  e^{\gamma D} Y(a,z)e^{-\gamma D}=i_{z;\gamma}Y(a,z+\gamma).
\end{equation}
Infinitesimally this gives another familiar formula: by
differentiating with respect to $\gamma$ we obtain
\begin{equation}
  \label{eq:InfTransCovYclass}
  [D,Y(a,z)]=\partial_zY(a,z).
\end{equation}
 Bicharacters
satisfying condition \eqref{eq:Ginvr} give rise to quantum vertex operator
algebras in the sense of Etingof-Kazhdan, \cite{MR2002i:17022}.  
In case the bicharacter satisfies \eqref{eq:Ginvr} and is also symmetric:
\[
r_{z_1,z_2}^\tau=r_{z_1,z_2}, 
\]
we obtain vertex operators of a vertex algebra as is usually
defined (see \cite{MR996026}, \cite{MR1651389}. This is a
special case of a more general result of Borcherds, see
\cite{MR1865087}, Theorem 4.2.

The condition \eqref{eq:Ginvr} is not satisfied in the case we are
interested in, see  section \ref{sec:mainexample}.

\section{Bicharacter Expansions and $S$-commutator}
\label{sec:BicharExpScom}

We continue to assume that $V$ has an $H_D$-quantum vertex algebra
structure via a bicharacter $r_{z_1,z_2}$, see Theorem
\ref{thm:h_d-quantum-vertex-from-bichar}. In this section we show how
an expansion of the bicharacter leads to a closed formula for the
$S$-commutator of fields.

Consider the vectorspace $V\otimes W(z)$, where $W(z)$ is some
space of functions (or power series) in $z$. Then we get an action of $H_D$
on this vector space by using the coproduct:
\[
D^{(k)} (a\otimes f(z))=\sum_{p+q=k} D^{(p)}a\otimes\partial_z^{(q)}f(z).
\]
\begin{thm}
  \label{thm:Scomandbichar}
  Let $a,b\in V$ and suppose that
\[\delta(r_{z_1,z_2}(a\otimes b))=\sum_{k\ge 0} d_k(a\otimes
b;t)\partial_{z_2}^{(k)}\delta(z_1,z_2),
\]
where $d_k(a\otimes b;t)\in \mathbb C[[z_1^{\pm 1},z_2^{\pm
  1}]][[t]]$. Then we have
\[
\comS{a(z_1),b(z_2)}=\sum_{k\ge 0} d_k(a^\prime\otimes
b^\prime;t)\sum_{p+q=k}Y\left([D^{(p)}a^{\prime\prime}]
  b^{\prime\prime},z_2\right)\partial^{(q)}_{z_2}\delta(z_1,z_2).
\]
\end{thm}

\begin{proof}
  The RHS of the $S$-commutator of the fields of $a$ and $b$ acting on
  $c$is
  \begin{align*}
    e^{z_1D}&a^\prime (e^{z_2D}b^\prime)c^\prime
    \delta\left(r_{z_1,z_2}({a^{\prime\prime}}^\prime\otimes
      {b^{\prime\prime}}^\prime)\right)
    r_{z_1,0}({a^{\prime\prime}}^{\prime\prime}\otimes
    {c^{\prime\prime}}^\prime)
    r_{z_2,0}({b^{\prime\prime}}^{\prime\prime}\otimes
    {c^{\prime\prime}}^{\prime\prime})=\\
    &= (e^{z_2D}b^\prime)c^\prime
    \sum_{k\ge0}d_k({a^{\prime\prime}}^\prime\otimes
    {b^{\prime\prime}}^\prime;t)
    \partial_{z_2}^{(k)}\left([e^{z_2D}a^\prime] r_{z_2,0}
      ({a^{\prime\prime}}^{\prime\prime}\otimes
      {c^{\prime\prime}}^\prime)\delta(z_1,z_2)\right)\times\\
    &\qquad\qquad\qquad \qquad\qquad\qquad \times r_{z_2,0}
    ({b^{\prime\prime}}^{\prime\prime}\otimes
    {c^{\prime\prime}}^{\prime\prime})=\\
    &= (e^{z_2D}b^\prime)c^\prime
    \sum_{k\ge0}d_k({a^{\prime\prime}}^\prime\otimes
    {b^{\prime\prime}}^\prime;t) \sum_{p+q+r=k}
    e^{z_2D}\left(D^{(p)}a^\prime\right)\times \\
    &\qquad\qquad \qquad\qquad\qquad \times
    r_{z_2,0}({D^{(q)}a^{\prime\prime}}^{\prime\prime}\otimes
    {c^{\prime\prime}}^\prime)r_{z_2,0}
    ({b^{\prime\prime}}^{\prime\prime}\otimes
    {c^{\prime\prime}}^{\prime\prime})\partial_{z_2}^{(r)}\delta(z_1,z_2)=
    \\
    &=\sum_{k\ge0}\sum_{p+q+r=k} d_k(a^\prime\otimes b^\prime;t)[
    e^{z_2D}\left(D^{(p)} {a^{\prime\prime}}^\prime\right)
    {b^{\prime\prime}}^{\prime} c^\prime] r_{z_2,0} ([D^{(q)}
    {a^{\prime\prime}}^{\prime\prime}]
    {b^{\prime\prime}}^\prime\otimes
    c^{\prime\prime})\times\\
    &\qquad\qquad\qquad \qquad\qquad\qquad \times \partial_{z_2}^{(r)}\delta(z_1,z_2)=\\
    &=\sum_{k\ge0}\sum_{p+q=k} d_k(a^\prime\otimes b^\prime;t) Y(
    (D^{(p)}a^{\prime\prime})b^{\prime\prime},z_2)c
    \partial_{z_2}^{(q)}\delta(z_1,z_2).
  \end{align*}   
\end{proof}

\section{The Main Example}
\label{sec:mainexample}

For the rest of the paper we will study a particular example of an
$H_D$-quantum vertex algebra $V$ obtained from a bicharacter as in Theorem
\ref{thm:h_d-quantum-vertex-from-bichar}. As an vector space $V$ is
the underlying space of the lattice vertex algebra based on the rank 1
lattice $\mathbb{Z}$ with pairing $(m,n)\mapsto mn$, cf.,
\cite{MR1651389}, section 5.4.

To define a bicharacter on $V$ we need an $H_D$-bialgebra
structure. As  $H_D$-bialgebra $V$ is generated by group-like elements
$e^\alpha,e^{-\alpha}$, so that
\[
\Delta(e^{m\alpha})=e^{m\alpha}\otimes e^{m},\quad
\epsilon(e^{m\alpha})=1,\quad m\in\mathbb Z.
\]
If we write $h=(De^\alpha)e^{-\alpha}$ then $h$ is primitive: we have
$\Delta(h)=h\otimes 1+1\otimes h$, $\epsilon(h)=0$. Then
\[
V=\bigoplus_{m\in \mathbb Z} V_m,\quad V_m= k[ D^{n}h]_{n\ge0}\otimes e^{m\alpha}.
\]
In fact $V$ is a Hopf algebra, with antipode $S\colon e^\alpha\mapsto
e^{-\alpha}$.  We define in this case a bicharacter on $V$ by putting
on generators
\begin{equation}
    r_{z_1,z_2}(e^{m\alpha}\otimes
    e^{n\alpha})=\sigma^{mn}, \quad \sigma=\frac{z_1-z_2}{1-tz_2/z_1},
\label{eq:DefBiChar}
\end{equation}
and extend to all of $V$ by using the properties of bicharacters, see
\cite{MR1865087} for details. Here (and below) we will expand any
rational expression in $t$ in \emph{positive} powers of $t$.  Note
that $r_{z_1,z_2}$ satisfies the VO assumption of Definition
\ref{defn:VO-assumption}. So by Theorem
\ref{thm:h_d-quantum-vertex-from-bichar} $V$ has an $H_D$-quantum
vertex algebra structure.

The bicharacter $r_{z_1,z_2}$ of this example is implicit in the paper
by Jing, \cite{MR1112626}.  By putting $t=0$ we obtain a bicharacter
$r^0_{z_1,z_2}$ which is implicit in the usual construction of a
lattice vertex algebra from the lattice $\mathbb Z$ with pairing
$(m,n)\mapsto mn$.

We will collect for later reference some values of this bicharacter
and of its associated braiding and translation bicharacters. First a
simple lemma.

\begin{lem}\label{lem:rhohexpmalpha_hh}
  For any bicharacter $\rho_{z_1,z_2}$ on $V$ we have, if
  $\rho_{z_1,z_2}(e^{m\alpha}\otimes e^{n\alpha})=\rho^{mn}$,
\[
\rho_{z_1,z_2}(h\otimes e^{m\alpha})=m\partial_{z_1}\ln(\rho), \quad
\rho_{z_1,z_2}(h\otimes h)=\partial_{z_2}\partial_{z_1}\ln(\rho).\]
\end{lem}

\begin{lem}\label{lem:rhexpmalpha_hh}
\begin{align*}
  r_{z_1,z_2}(h\otimes
  e^{m\alpha})&=m\left(\frac1{z_1-z_2}-\frac{tz_2/z_1}
    {z_1-tz_2}\right),\\
  r_{z_1,z_2}(h\otimes h)&=\frac1{(z_1-z_2)^2}-\frac
  t{(z_1-tz_2)^2}.
\end{align*}
\end{lem}
The bicharacter $r_{z_1,z_2}$ is invertible ($V$ being a Hopf
algebra), with inverse on generators given by
\[
r\inv_{z_1,z_2}(e^{m\alpha}\otimes e^{n\alpha})=\sigma^{-mn}.
\]
\begin{lem}\label{lem:braidinghexpmalpha_hh}
The 
braiding bicharacter $R_{z_1,z_2}$ of $r_{z_1,z_2}$ is given on the
generators by
\begin{equation}
  \label{eq:BradingGen}
  R_{z_1,z_2}(e^{m\alpha}\otimes e^{n\alpha})=\Sigma^{mn},\quad
  \Sigma=\Sigma_{z_1,z_2}=-\frac{1-tz_2/z_1}{1-tz_1/z_2},
\end{equation}
and we have
\begin{equation}
  \label{eq:braidinghexphh}
  R_{z_1,z_2}(h\otimes
  e^{m\alpha})=m\left(\frac{tz_2/z_1}{z_1-tz_2}+
\frac{t}{z_2-tz_1}\right),
\end{equation}
and 
\begin{equation}
  \label{eq:Bradinghh}
    R_{z_1,z_2}(h\otimes
    h)=\frac{t}{(z_1-tz_2)^2}-\frac t{(z_2-tz_1)^2}.
\end{equation}
\end{lem}

\begin{lem}\label{lem:Translhexpmalpha_hh}
  The translation bicharacter $R^{\gamma}$ of $r_{z_1,z_2}$ is given
  on generators by
\begin{equation}
  \label{eq:TranslGen}
  R^\gamma_{z_1,z_2}(e^{m\alpha}\otimes e^{n\alpha})=\Pi^{mn},\quad
  \Pi=\Pi_{z_1,z_2}=\frac{1-tz_2/z_1}{1-t\frac{z_2+\gamma}{z_1+\gamma}},
\end{equation}
and we have
\begin{equation}
  \label{eq:transhexphh}
  R^\gamma_{z_1,z_2}(h\otimes
  e^{m\alpha})=\frac{mtz_2/z_1}{z_1-tz_2}-
  \frac{mt({z_2+\gamma })/({z_1+\gamma })}{(z_1+\gamma )- t(z_2+\gamma )},
\end{equation}
and 
\begin{equation}
  \label{eq:Translhh}
    R^\gamma_{z_1,z_2}(h\otimes    h)
    =\frac t{(z_1-tz_2)^2}-\frac t{((z_1+\gamma )-t(z_2+\gamma ))^2}.
\end{equation}
\end{lem}

We will calculate some $(n)$-products of states and of fields in $V$ to
illustrate what is involved.

First note that $r_{z,0}(e^\alpha\otimes e^{-\alpha})=\frac 1z$. 
This implies that
\[
Y(e^\alpha,z)e^{-\alpha}=(e^{zD}e^\alpha)e^{-\alpha}r_{z,0}(e^\alpha\otimes
e^{-\alpha})=\frac1z+h+\mathcal{O} (z),
\]
so that we have the following products of states.
\begin{equation}
e^\alpha_{(-1)}e^{-\alpha}=h,\quad e^\alpha_{(0)}e^{-\alpha}=1, \quad
e^\alpha_{(k)}e^{-\alpha}=0, k>0.\label{eq:expexpproducts}
\end{equation}
Note that this are the same $(n)$-products as for the lattice vertex
algebra corresponding to the bicharacter $r^0_{z_1,z_2}$ (obtained by
putting $t=0$).

Next we want to use Corollary \ref{cor:STFieldhom} to calculate
$(n)$-products of fields. We have by Lemma
\ref{lem:Translhexpmalpha_hh} $R^{z_2}_{z_3,0}(e^\alpha\otimes
e^{-\alpha})=1-t\frac{z_3}{z_2+z_3}$, so that
\begin{align*}
  i_{z_2;z_3}S^{(z_2)}_{z_3,0}(e^\alpha\otimes
  e^{-\alpha})&=e^\alpha\otimes
  e^{-\alpha}i_{z_2;z_3}(1-t\frac{z_3}{z_2+z_3}) =\\
&=e^\alpha\otimes
  e^{-\alpha}(1+t\sum_{k=1}^\infty (-1)^k (\frac{z_3}{z_2})^k).
\end{align*}
Hence by Corollary \ref{cor:STFieldhom} and \eqref{eq:expexpproducts}
\begin{align*}
  e^\alpha(z)_{(-1)}e^{-\alpha}(z)&= Y(e^\alpha_{(-1)}e^{-\alpha},z)
  -Y(e^\alpha_{(0)}e^{-\alpha},z)\frac tz\\
  &=h(z)- \frac tz.
\end{align*}
Now, see Section \ref{sec:NormOrdProd}, $e^\alpha(z)_{(-1)}e^{-\alpha}(z)=\nopS{e^\alpha(z)e^{-\alpha}(z)}$,
and this normal ordered product of fields is \emph{not} a vertex
operator $Y(a,z)$ for any $a\in V$, since the action of
$\nopS{e^\alpha(z)e^{-\alpha}(z)}$ on the vacuum is not regular in
$z$, contradicting the vacuum axiom \eqref{eq:vacuumX}. This is in
contrast to the situation in the usual vertex algebras.

\section{$S$-Commutators and Commutators}
\label{sec:ScomandCom}

In this section we calculate some $S$-commutators of fields by
expanding the bicharacter in our main example and express this in
terms of commutators,  using Theorem \ref{thm:Scomandbichar}. 

We have
\[
\delta\Big(r_{z_1,z_2}(e^{m\alpha}\otimes e^{n\alpha})\Big)=
\begin{cases}
  0& mn\ge0,\\
(1-t{z_2}/{z_1})^{k+1}\partial_{z_2}^{(k)}\delta(z_1,z_2) &mn=-k-1<0,
\end{cases}
\]
which follows from the definition \eqref{eq:DefBiChar}.  Then
\[
\comS{e^{m\alpha}(z_1),e^{n\alpha}(z_2)}=
  \begin{cases}
    0& mn\ge0\\
    (1-\frac{t{z_2}}{z_1})^{k+1}\sum Y(v^p_{m,n},z_2)
    \partial_{z_2}^{(q)}\delta(z_1,z_2) &mn=-k-1<0
      \end{cases}
\]
where $v^p_{m,n}=D^{(p)}\left(e^{m\alpha}\right)e^{n\alpha}\in V$ and
the sum is over all $p,q\ge0$ such that $p+q=k$. In particular 
\[
\comS{e^{\alpha}(z_1),e^{-\alpha}(z_2)}=(1-t\frac
{z_2}{z_1})\delta(z_1,z_2)=(1-t)\delta(z_1,z_2).
\]
So
\[
e^{\alpha}(z)_{(0)}e^{-\alpha}(z)=1-t,\quad
e^{\alpha}(z)_{(k)}e^{-\alpha}(z)=0, k>0.
\]
In the same way 
\[
\delta\Big(r_{z_1,z_2}(h\otimes e^{m\alpha})\Big)=m\delta(z_1,z_2),
\]
which follows from Lemma \ref{lem:rhexpmalpha_hh}, see also
(\ref{eq:deltaz-tw})). Hence
\begin{equation}
\comS{h(z_1),e^{m\alpha}(z_2)}=me^{m\alpha}(z_2)\delta(z_1,z_2).\label{eq:comhexp}
\end{equation}
Finally, 
using Lemma \ref{lem:rhexpmalpha_hh} again, we find
\[
\delta\Big(r_{z_1,z_2}(h\otimes h)\Big)=\partial_{z_2}\delta(z_1,z_2),
\]
so that, 
\begin{equation}
\comS{h(z_1),h(z_2)}=\partial_{z_2}\delta(z_1,z_2).\label{eq:comhh}
\end{equation}
It is sometimes useful to express the $S$-commutators of fields in
terms of the usual commutators. We give some examples.

We have by definition of the $S$-commutator
 \begin{align*}
   \comS{h(z_1),e^{m\alpha}(z_2)}&=h(z_1)e^{m\alpha}(z_2)-e^{m\alpha}(z_2)h^\prime(z_1)
 R_{z_2,z_1}(e^{m\alpha}\otimes h^{\prime\prime})=\\
 &=[h(z_1),e^{m\alpha}(z_2)]-e^{m\alpha}(z_2)R_{z_2,z_1}(e^{m\alpha}\otimes
 h)=\\
 &=[h(z_1),e^{m\alpha}(z_2)]-e^{m\alpha}(z_2)m\partial_{z_1}\ln(\Sigma_{z_2,z_1}),
 \end{align*}
 where $\Sigma$ is defined in Lemma \ref{lem:braidinghexpmalpha_hh}.
 Combining this with \eqref{eq:comhexp} gives
\begin{equation}
  \label{eq:Comh_emalpha}
  [h(z_1),e^{m\alpha}(z_2)]=me^{m\alpha}(z_2)\big(\delta(z_1,z_2)+\partial_{z_1}\ln(\Sigma_{z_2,z_1})\big).
\end{equation}
Now 
\[
\Res_{z_1}\Big( z_1^n\partial_{z_1}\ln(\Sigma_{z_2,z_1})\Big)=
  \begin{cases}
    0&n=0,\\
-t^{\abs{n}}z_2^n& n\ne0.
  \end{cases}
\]
Hence
\begin{equation}
  [h_{(n)},e^{m\alpha}(z_2)]=
  \begin{cases}
    me^{m\alpha}(z_2)&n=0,\\
mz_2^n(1-t^{\abs{n}})e^{m\alpha}(z_2)&n\ne0.
  \end{cases}\label{eq:comhcompealpha}
  \end{equation}
Similarly,
\begin{align*}
  \comS{h(z_1),h(z_2)}&=[h(z_1),h(z_2)]-R_{z_2,z_1}(h\otimes
  h)=\\
&=[h(z_1),h(z_2)]-\big(\frac t{(tz_1-z_2)^2}-\frac
  t{(tz_2-z_1)^2}\big),
\end{align*}
Note that here we see that the ordinary commutator of $h(z)$ with
itself is not killed by any power of $z_1-z_2$, whereas the
$S$-commutator is killed by $(z_1-z_2)^2$, see \eqref{eq:comhh}.  
 
By \eqref{eq:comhh}
\[
[h(z_1),h(z_2)]=\partial_{z_2}\delta(z_1,z_2)+R_{z_2,z_1}(h\otimes h).
\]
Now
\[
\Res_{z_1}\Big( z_1^n R_{z_2,z_1}(h\otimes h)\Big)=
-nt^{\abs{n}}z_2^{n-1},
\]
and we have
\begin{equation}
  \label{eq:hnhz}
  [h_{(m)},h(z_2)]=mz_2^{m-1}(1-t^{\abs m})
\end{equation}
and 
\begin{equation}
  \label{eq:hmhn}
  [h_{(m)},h_{(n)}]=m(1-t^{\abs m})\delta_{m+n,0}.
\end{equation}
We see therefore that the coefficients of $h(z)$ generate a
\emph{deformed Heisenberg algebra} $\mathcal{H}_t$. As a Lie algebra
$\mathcal{H}_t$ is isomorphic to the usual Heisenberg Lie algebra
$\mathcal{H}=\mathcal{H}_{t=0}$. In particular the representation
theory of $\mathcal{H}_t$ is the same as in the undeformed case.
We have a decomposition
\[
V=\oplus_{m\in \mathbb Z}V_m,\quad V_m=k[D^nh]e^{m\alpha},
\]
where each $V_m$ is an irreducible $\mathcal{H}_t$-module, with action
given by
\[
h_{(m)}=
\begin{cases}
  \text{multiplication by } D^kh/k!& m=-k-1<0,\\
\partial_\alpha& m=0,\\
m(1-t^m)\frac{\partial}{\partial h_{(-m-1)}}& m>0.
\end{cases}
\]
The case $m=0$ follows from Cor \ref{cor:STFieldhom} and \eqref{eq:comhexp}.

\section{Braided Bosonization}
\label{sec:braidedboson}
Define
\[
\Gamma_+(z)=\exp\left(\sum_{n>0}h_{(-n)}z^n/n\right),\quad\Gamma_-(z)=
\exp\left(-\sum_{n>0}h_{(n)}z^{-n}/n\right).
\]
By \eqref{eq:hmhn} we have for $m\ne0$
\[
[h_{(\pm m)}, \Gamma_\pm(z)]=\pm z^{\pm m}(1-t^{\abs m})\Gamma_\pm(z),
\quad [h_{(\mp m)}, \Gamma_\pm(z)]=0.
\]
Then we see that
\[
\Sigma_n(z)=\Gamma_+^{-n}(z)e^{n\alpha}(z)\Gamma_-^{n\alpha}(z)e^{n\alpha}
\]
commutes with the deformed boson:
\[
[h(z_1),\Sigma_n(z_2)]=0,
\]
and by the usual arguments using the representation theory of the
deformed Heisenberg algebra (see e.g., \cite{MR1651389}) one  finds the
bosonization formula
\[
e^{n\alpha}(z)=\Gamma_+^n(z)\Gamma_-^{-n}(z)e^{n\alpha}z^{n\partial_\alpha}.
\]
This formula (for $n=\pm1$) can by found in Jing's paper, \cite{MR1112626}, with a
slightly different notation.

\section{Hall-Littlewood Polynomials}
\label{sec:hall-littl-polyn}

In this section we recall the Macdonald definition of Hall-Littlewood
symmetric polynomials (\cite{MR1354144}). Also we
explain how the bosonized vertex operators described in the previous
secton (as considered by N. Jing, \cite{MR1112626}), serve as
generating functions for the Hall-Littlewood poynomials.

Denote by $\Lambda $ the ring of symmetric functions over
$\mathbf{C}[[t]]$ in countably many independent variables $x_i,\ i \ge
0 $.

Let $\lambda$ be a partition, $\lambda =(\lambda _1, \lambda _2,
\dots,\lambda _k, \dots), \ \lambda _1 \ge \lambda _2 \ge \dots \ge
\lambda _k \ge\dots$.  Let $\abs{\lambda }=\lambda _1 +\lambda _2
+\dots +\lambda _k +\dots $.

Denote $z_{\lambda}=\prod _{i\ge 0}i^{m_i}.m_i!$, where
$m_i=m_i(\lambda)$ is the number of parts of $\lambda$ equal to $i$.

We call a family $(a_{\lambda})$ of elements in a ring indexed by
partitions \emph{multiplicative} if $a_{\lambda}=\prod a_{\lambda
  _i}$.

For any partition $\alpha $ we use the vector notation $x^{\alpha }$
for $x_1^{\alpha _1}x_2^{\alpha _2}\dots x_k^{\alpha _k}\dots $.  We
will use the basis $(m_{\lambda})$ of monomial symmetric functions:
\begin{equation}
m_{\lambda}=\sum_{\alpha } x^{\alpha },
\end{equation}
where the sum is over distinct permutations of $\lambda $, as well as
the multiplicative basis generated by the power sums $p_n=\sum_{i\ge
  0}x_i^n, \ p_0=1$.

Define a scalar product $\langle \ ,\ \rangle _{t}$ on $\Lambda
_{\mathbf{F}}$ by putting for the power functions
\begin{equation*}
\label{eq:scalarpr}
\langle p_{\lambda}, p_{\mu} \rangle _{t} =\delta _{\lambda \mu} z_{\lambda}v_{\lambda},
\end{equation*}
for any partitions $\lambda ,\mu$, where the mulitplicative family
$v_{\lambda}$ is defined by $v_n=\frac{1}{1-t^n}$.  Define a set of
symmetric functions $\{H_{\lambda}\}$ indexed by partitions by the
following two (over-determining) conditions:
\begin{align*}
  &\langle H_{\lambda}, H_{\mu} \rangle _{t}=0 \ \ \ \text{for} \ \ \lambda \neq \mu, \\
  &H_{\lambda}=m_{\lambda}+\sum _{\mu < \lambda}u_{\lambda
    \mu}m_{\mu}, \ \ u_{\lambda \mu}\in \mathbf{C}[[t]].
\end{align*}
Here $\mu < \lambda$ is with respect to the usual partial order on
partitions.

It is proved in (\cite{MR1354144}) that such symmetric functions
$\{H_{\lambda}\}$ exist.  Denote also by $Q_{\lambda}$ the dual of
$H_{\lambda}$, i.e., $\langle H_{\lambda},Q_{\mu}\rangle _{t}=\delta
_{\lambda , \mu }$. Note that when $t=0$ both the $H_{\lambda}$ and
the $Q_{\lambda}$ reduce to the Schur polynomials (Schur poynomials
are self dual).

We can view the $n$-th power symmetric function $p_n$ as an operator
acting on $\Lambda _\mathbf{C}[[t]]$ by multiplication. Define also
for given multiplicative family $(v_{\lambda })$ the operators
$p^{\perp}_{n}$ by requiring
\begin{displaymath}
\langle p^{\perp}_{n}f,g\rangle _{t}= \langle f,p_n g\rangle _{t},
\end{displaymath}
for any $f,g \in \Lambda _\mathbf{C}[[t]]$.

\begin{lem}
\label{lem:heis}
The operators $\{ h_{(n)} | n\in \mathbf{Z}\}$ given by
\mbox{$h_{(n)}=-(1-t^{n})p^{\perp}_n$,} \\
\mbox{$h_{(-n)}=(1-t^{n})p_n$} \ for $n\in \mathbf{N}$, \ $h_{(0)}=0$
generate a representation of the deformed Heisenberg algebra
$\mathcal{H}_t$ on $\Lambda _\mathbf{C}[[t]]$, i.e.,
\begin{equation}
 [h_{(m)},h_{(n)}]=m(1-t^{\abs m})\delta_{m+n,0}.
\end{equation}
\end{lem}
The proof is based on the undeformed case ($v_n=1$), which can be
found in \cite{MR1354144}.

From the fact that the power symmetric functions form a basis of
$\Lambda _\mathbf{C}[[t]]$, it follows that $\Lambda _\mathbf{C}[[t]]$
is a highest weight module for $\mathcal{H}_t$, and is thus an
irreducible $\mathcal{H}_t$ module. Therefore we have that $\Lambda
_\mathbf{C}[[t]]$ is isomorphic as a module and as an algebra to
$V_0$ ($V_0$ was defined in Section \ref{sec:mainexample}). Thus we can identify $(1-t^{n})p_n$ with $D^{(n-1)}h \ (n>0)$.

The following theorem (\cite{MR1112626}) explains the connection
between the Hall\--Little\-wood symmetric functions and the vertex
operators considered in the previous section:
\begin{thm}
  Let $\tilde{m}$ is a partition of length $l$, $\tilde{m}=(m_1, m_2,
  \dots , m_l, 0, \dots)$, and let $\rho$ be the partition defined by
  $\rho=(l, l-1, \dots ,1, 0, \dots )$). The constant term of
  \mbox{$Y(D^{(m_1)}e^{\alpha }, z_1)Y(D^{(m_2)}e^{\alpha }, z_2)\dots
    Y(D^{(m_1)}e^{\alpha }, z_l)1$} is
  $Q_{\tilde{m}-\rho}e^{l\alpha}$, where $Q_{\tilde{m}-\rho}$ is the
  dual Hall-Littlewood polynomial corresponding to the partition
  ${\tilde{m}-\rho}$.
\end{thm}
The proof is straightforward modification of the main theorem in
\cite{MR1112626} using the properties of the vertex operators.

Thus the vertex operators $Y(D^{(m)}e^{\alpha }, z)$ (as described in
Section \ref{sec:mainexample}) and the coefficients of their
products are very important in the theory of the Hall-Littlewood
poynomials.  This makes them an important example of
quantum vertex operators, and they are the main motivation for our
definition of $H_D$-quantum vertex algebras. The previous definitions
of quantum vertex algebras were not general enough to incorporate the
Hall-Littlewood vertex operators.

\appendix

\section{Braided Algebras with symmetry}
\label{sec:BraidAlgSym}
 
\subsection{Introduction}
\label{subsec:AppIntr}

To motivate the rather complicated definition of an $H_D$-quan\-tum
vertex algebra in Section \ref{sec:H_DQuantumVertex} we discuss in
this Appendix braided algebras (with symmetry). The idea is that a
vertex algebra has a singular multiplication, and that it is good to
understand the nonsingular case first.

\subsection{Commutative Associative Algebras}
\label{ubsec:AppComAs}

As a preliminary, note that an efficient way to describe commutative
associative unital algebras is as follows. Let $M$ be a vector space
and $1\in M$ a distinguished element, and let
\[
m\colon M^{\otimes 2}\to M
\]
be a multiplication for which $1$ is the unit:
\begin{equation}
  \label{eq:unitalalg}
  m(a\otimes 1)=m(1\otimes a)=a, \quad a\in M.
\end{equation}
We need some notation. If $a$ is a linear map on $M^{\otimes 2}$, then
$a^{23}$ is the operator on $M^{\otimes 3}$ acting on the $2$nd and
$3$rd factor (so $a^{23}=1\otimes a$). The other superscripts have a
similar meaning. Let $\tau\colon M^{\otimes 2}\to M^{\otimes 2}$ be
the flip $a\otimes b\mapsto b\otimes a$. Let $m_3=m(1\otimes m)\colon
M^{\otimes 3}\to M$. Then we impose
\begin{equation}
  \label{eq:ComAx}
  m_3= m_3\tau^{12},\tag{\textbf{Commutativity/Associativity Axiom}}
\end{equation}
In other words, writing $m_3(a\otimes b\otimes c)=a(bc)$, we require
$a(bc)=b(ac)$. Then one easily checks that $(M,m,1)$ is in fact
commutative ($m=m\tau$) and associative ($m(1\otimes
m)=m(m\otimes1)$).

\subsection{Braided Algebras}
\label{subsec:Braidedalg}
We are next interested in non commutative algebras where the
noncommutativity is controlled by a braiding map.

\begin{defn}\label{defn:Braidedalg}
  A \emph{braided algebra} is a unital algebra $(M,m,1)$ with a
  braiding $S\colon M^{\otimes
  2}\to M^{\otimes 2}$ such that
\begin{enumerate}
\item \upshape{(\textbf{Vacuum Axiom})} $S(a\otimes 1)=a\otimes 1$,
  and $S(1\otimes a)=1\otimes a$.
\item \upshape{(\textbf{Braiding Axiom})} $m_3 S^{12}=m_3 \tau^{12}$.
\item \upshape{(\textbf{Unitarity Axiom})} $ S\circ \tau\circ S\circ \tau=1_{M^{\otimes 2}}$.
\item \upshape{(\textbf{Yang-Baxter Axiom})} $ S^{12}S^{13}S^{23}=S^{23}S^{13}S^{12}$.
\item \upshape{(\textbf{Compatibility with Multiplication Axiom})}
  $Sm^{12}=m^{12}S^{23}S^{13}$ and $Sm^{23}=m^{23}S^{12}S^{13}$.
\end{enumerate}
\end{defn}

The Compatibility with Multiplication Axiom allows us to express the
braiding involving a product in terms of a product of the braidings of
the factors. Also, together with the braiding axiom it gives associativity,
as we now proceed to show.

\begin{lem}[\textbf{Braided Commutativity}] \label{lemA:braid-com}
\[m S= m \tau.\]  
\end{lem}
\begin{proof}
Apply the Braiding Axiom to $a\otimes b\otimes 1$, using $m_3(a\otimes
b\otimes 1)=m(a\otimes b)$.
\end{proof}

\begin{thm}[\textbf{Associativity}] \label{thm:assoc} 
A braided algebra is associative:
\[
m (1\otimes m)=m(m\otimes 1).
\]
\end{thm}

\begin{proof}
 Let $A=a\otimes b \otimes\ c$. Then
  \begin{align*}
    m m^{12}(A)&=m S\tau m^{12}(A)=m S m^{23}(c\otimes a\otimes b)=\\
    &=m m^{23}S^{12}S^{13}(c\otimes a\otimes b)=m m^{23}\tau^{12}S^{13}(c\otimes a\otimes b)=\\
    &=m m^{23}S^{23}(a\otimes c\otimes b)=m m^{23}\tau^{23}(a\otimes c\otimes b)=\\
    &=m m^{23}(A)
  \end{align*}

\end{proof}

We used the Compatibility with Multiplication Axiom to derive
associativity. If we don't impose this axiom, we can only derive
\emph{braided associativity}, (also called \emph{quasi-associativity}
cf. \cite{MR2002i:17022}):
\[
  m m^{23}S^{23}S^{13}=m S m^{12},\quad m m^{12}S^{12}S^{13}=m S
  m^{23}.
\]
We have not yet used the unitarity and Yang-Baxter axioms. They are
used to describe the behaviour under permutations of the arguments
of the $n$-fold multiplication $m_n\colon M^n\to M$ (defined
recursively by $m_n=m(1\otimes m_{n-1})$) as we now proceed to
explain.
\begin{lem}\label{lem:simpletransponm}
  \[
m_n \tau^{ii+1}=m_n S^{ii+1}.
\]
\end{lem}
\begin{proof}
  We can use associativity to write
\[
m_n=m_3\circ(m_{i-1}\otimes m_2 \otimes m_{n-i-1}).
\]
The Lemma follows from Braided Commutativity, Lemma \ref{lemA:braid-com} .
\end{proof}

\begin{remark}
Note that if $i, j$ are not adjacent, then it is in general not true that
 the transposition $\tau^{ij}$ does  act on $m_n$ by multiplication
 by $S^{ij}$.

 For instance, a simple example of a non trivial braided algebra is a
 super commutative algebra $M=M_{\bar 0}\oplus M_{\bar 1}$. The
 braiding is given (for homogeneous elements)  by $S(a\ten
 b)=(-1)^{\abs a\abs b}a\ten b$. 
It is then clear that the braiding corresponding to the permutation
$\tau^{13}\colon a\ten b\ten c\mapsto c\ten b\ten a$ is given by
\[
S^{\tau^{13}}(a\ten b\ten c)=S^{12} S^{13}S^{23}(a\ten b\ten
c)=(-1)^{\abs a\abs b}(-1)^{\abs a\abs c}(-1)^{\abs b\abs c}(a\ten
b\ten c),
\]
whereas 
\[
S^{13}(a\ten b\ten c)=(-1)^{\abs a\abs c}a\ten b\ten c.
\]\qed
\end{remark}

One knows that the symmetric group $\mathcal{S}_n$ is generated by the
simple transpositions $w_i=(ii+1)$, $i=1, 2, \dots, n-1$, see Section
\ref{sec:braidingmaps}. Then define a map $S\colon \mathcal{S}_n\to
\GL(M^{\ten n})$ by
\[
S(w_i)=1^{i-1}\ten S\tau\ten 1^{n-i-1},
\]
and extend this as an anti-homomorphism:
\[
S(f)=S(w_{i_k})S(w_{i_{k-1}})\dots S(w_{i_1}),
\]
in $f=w_{i_1}w_{i_2}\dots w_{i_k}\in\mathcal{S}_n$.
Then the unitarity and the Yang-Baxter axioms and Lemma
\ref{lem:simpletransponm} imply
\begin{thm}
  The braiding map $S^{f}:M^{\otimes n}\to M^{\otimes n}$ is
  independent of the representation of $\sigma$ in terms of simple
  reflections. Furthermore
\[
m_n \circ S(f)= m_n .
\]
for all $f\in \mathcal{S}_n$.
\end{thm}
This concludes our discussion of braided algebras \emph{an sich}.

\subsection{Braided Algebras with Symmetry}
\label{subsec:BrRiSymm}

We now assume that we have additionally an action of a group $G$ on
the braided algebra $M$. If $g\in G$ we write $\Delta(g)=g\otimes g\in
G\otimes G$ for the coproduct of $g$.

\begin{defn}\label{defn:BraidedGalg}
  Let $(M,m,1, S)$ be a braided algebra, with a $G$-action on $M$. We
  call this a braided $G$-algebra in case for each $g\in G$ there is a
  map
\[
S^g:M^{\otimes 2}\to M^{\otimes 2},
\]
such that
\begin{itemize}
\item \upshape{(\textbf{Vacuum Axiom})} $S^g(a\otimes 1)=a\otimes 1$,
  and $S^g(1\otimes a)=1\otimes a$.
\item \upshape{(\textbf{$G$-Symmetry})} $g m S^g=m \Delta(g)$.
\item \upshape{(\textbf{Multiplicativity})} $S^{gh}=S^h\circ
  \Delta(h\inv)\circ S^g \circ\Delta(h)$
\item \upshape{(\textbf{$G$-Yang-Baxter})} $ S^{g,12}S^{g,13}S^{g,23}=S^{g,23}S^{g,13}S^{g,12}$.
\item \upshape{(\textbf{Compatibility with Multiplication Axiom})}
  $S^gm^{12}=m^{12}S^{g,23}S^{g,13}$ and  $S^gm^{23}=m^{23}S^{g,12}S^{g,13}$.
\end{itemize}
\end{defn}

Of course, the simplest case is were $S^g=1\otimes 1$ for all $g\in G$. Then
the multiplication intertwines the action of $G$ on $M^{\otimes 2}$
and $M$; usually $M$ is then called a module-algebra.

\begin{lem}
  Define $\Sigma_n,\tilde \Sigma_n\colon M^{n+1}\to M^{n+1}$ by
$\Sigma_n=S^{12}S^{13}\dots S^{1n+1}$, $\tilde\Sigma_n=S^{1n+1}\dots
S^{13} S^{12}$. Then we have compatibility with the higher multiplications:
\[
  S(1\otimes m_n)=(1\otimes m_n)\Sigma_n,\quad   S(m_n\otimes 1)=(
  m_n\otimes 1)\tilde\Sigma_n.
\]
\end{lem}
\begin{proof}
  For $n=2$ the Lemma is just the compatibility with multiplication
  axiom. Assume the Lemma is true for n=k-1. Then
  \begin{align*}
    S(1\otimes m_k)&= S(1\otimes m)(1\otimes1\otimes m_{k-1})=\\
                   &= (1\otimes m)S^{12}S^{13}(1\otimes 1\otimes
                   m_{k-1})=\\
                   &=(1\otimes m)(1\otimes1\otimes
                   m_{k-1})S^{12}\Sigma_{n-1}^{13\dots nn+1},
  \end{align*}
  Noting that $\Sigma_n=S^{12}\circ \Sigma_{n-1}^{13\dots nn+1}$ the
  first equation of the Lemma follows. The second one is proved
  similarly.
\end{proof}

Now define
\[
S^g_n=\Sigma_{n-1}\circ (1\otimes S^g_{n-1}).
\]
\begin{thm}
  We have $S^g_n=\tilde \Sigma_{n-1}\circ (S^g_{n-1}\otimes 1))$ and
\[
g m_n S_n^g=m_n \Delta_n(g).
\]
\end{thm}
\begin{proof}
  ??
\end{proof}

\subsection{Bicharacters}
\label{subsec:Bichar}
Let $M$ be a commutative and cocommutative Hopf algebra. A
bicharacters on $M$ is a linear map
\[
r\colon M^{\otimes 2}\to \mathbb{C},
\]
satisfying
\begin{itemize}
\item (\textbf{Vacuum}) $r(a\otimes 1)=r(1\otimes a)=\epsilon(a)$,
  $a\in M$.
\item (\textbf{Multiplication}) For all $a,b,c\in M$ we have
  $r(a\otimes bc)=\sum r(a^\prime\otimes b)r(a^{\prime\prime}\otimes
  c)$ and $r(ab\otimes c)=\sum r(a\otimes c^\prime)r(b\otimes
  c^{\prime\prime})$.
\end{itemize}
Here and below we use the notation $\Delta(a)=\sum a^\prime\otimes
a^{\prime\prime}$ for the coproduct for $a\in V$. Often we will also omit the
summation symbol, to unclutter the formulas.

We can multiply bicharacters: if $r,s$ are bicharacters and $a,b\in M$ then
\begin{equation}
  (r\ast s)(a\otimes b)=r(a^\prime\otimes b^\prime)
  s(a^{\prime\prime}\otimes b^{\prime\prime
  }).\label{eq:Appmultbichar}
\end{equation}
The unit bicharacter is 
\begin{equation}
\epsilon(a\otimes
b)=\epsilon(a)\epsilon(b).\label{eq:Appunitbichar}
\end{equation}
Since $M$ is a Hopf algebra it comes with an antipode, and all
bicharacters are invertible, with inverse given by
\[
r\inv(a\otimes b)=r(S(a)\otimes b).
\]
The set of bicharacters forms an Abelian group.

The transpose of a bicharacter is defined by
\[
r^\tau(a\otimes b)=r(b\otimes a).
\]
The transpose is an involution of the algebra of bicharacters:
\[
(r\ast s)^\tau =(r^\tau\ast s^\tau).
\]
If $r$ is an invertible bicharacter with inverse $r\inv$ we define another
bicharacter
\begin{equation}
  \label{eq:AppdefR}
R=r^\tau\ast r\inv,
\end{equation}
We will call $R$ the \emph{braiding bicharacter} associated to $r$. It
will control the braiding in the braided algebra we are going to
construct from $r$ below.  The braiding bicharacter $R$ is
\emph{unitary}:
\begin{equation}
  \label{eq:AppRunitary}
  R^\tau=R\inv.
\end{equation}
Also we have
\begin{equation}
  \label{eq:Appbraidingforr}
r\ast R=r^\tau.
\end{equation}
For any bicharacter $\rho$ on $M$ we define a map
\begin{equation}
S^{(\rho)}\colon M^{\otimes 2}\to M^{\otimes 2}, \quad a\otimes
b\mapsto a^\prime\otimes b^\prime \rho(a^{\prime\prime}\otimes
b^{\prime\prime}).\label{eq:DefSrho}
\end{equation}
\begin{lem}
  \begin{enumerate}
  \item If $\epsilon$ is the unit bicharacter on $M$, then
  $S^{(\epsilon)}=1_{M^{\otimes2}}$.
\item If $\rho,\sigma$ are bicharacters on $M$, then
  $S^{(\rho\ast\sigma)}=S^{(\rho)}\circ S^{(\sigma)}$.
\item If $\rho$ is a bicharacter, then $\tau\circ S^{(\rho)}\circ\tau=S^{(\rho^\tau)}$.
  \end{enumerate}
\end{lem}

\begin{lem}\label{lem:SrhoVacYangBaxter}
  For all $a\in M$ and bicharacters $\rho$ on $M$ we have
  \begin{enumerate}
  \item (\textbf{Vacuum}) $S^{(\rho)}(a\otimes 1)=a\otimes 1$ and
  $S^{(\rho)}(1\otimes a)=1\otimes a$.
\item (\textbf{Yang-Baxter}) $ S^{(\rho),12}S^{(\rho),13}S^{(\rho),
    23}=S^{(\rho),23}S^{(\rho),13}S^{(\rho),12}$.
  \end{enumerate}
\end{lem}
Now we fix a bicharacter $r$ on $M$, and define a twisting of the
multiplication $m$ on $M$:
\[
m_r=m\circ S^{(r)}\colon M^{\otimes 2}\to M.
\]
\begin{lem}\label{lem:rhocompatmr}
  For any bicharacter $\rho$ the map $S^{(\rho)}$ is compatible with
  the twisted multiplication $m_r$:
\[
S^{(\rho)}(m_r\otimes 1)=(m_r\otimes
1)S^{(\rho),23}S^{(\rho),13},\quad 
S^{(\rho)}(1\otimes m_r)=(1\otimes m_r)S^{(\rho),12}S^{(\rho),13}
\]
\end{lem}
\begin{proof}
  For $a,b,c\in M$ we have
  \begin{align*}
    S^{(\rho)}(m_r\otimes 1)(a\otimes b\otimes c)&=
    S^{(\rho)}(a^\prime b^\prime\otimes c)r(a^{\prime\prime}\otimes
    b^{\prime\prime})=\\
    &=(a^\prime b^\prime)^\prime\otimes c^\prime\rho((a^\prime
    b^\prime)^{\prime\prime}\otimes c^{\prime
      \prime})r(a^{\prime\prime}\otimes
    b^{\prime\prime})=\\
    &={a^\prime}^\prime {b^\prime}^\prime\otimes c^\prime
    \rho({a^\prime}^{\prime\prime}\otimes {c^{\prime\prime}}^\prime)
    \rho({b^\prime}^{\prime\prime}\otimes {c^{\prime
        \prime}}^{\prime\prime})
    r(a^{\prime\prime}\otimes  b^{\prime\prime})\\
 \intertext{Now by coassociativity and cocommutativity of $M$ we       have 
 ${a^\prime}^\prime\otimes {a^\prime}^{\prime\prime}\otimes
       a^{\prime\prime}={a^\prime}^\prime\otimes \otimes
       a^{\prime\prime}\otimes {a^\prime}^{\prime\prime}$, 
so that we get}
     &={a^\prime}^\prime {b^\prime}^\prime
    r({a^\prime}^{\prime\prime}\otimes {b^\prime}^{\prime\prime} )\otimes c^\prime
    \rho(a^{\prime\prime}\otimes {c^{\prime\prime}}^\prime)
    \rho( b^{\prime\prime}\otimes {c^{\prime
        \prime}}^{\prime\prime})=
    \\
    &=(m_r\otimes 1)S^{(\rho),23}S^{(\rho),13}(a\otimes b\otimes c).
\end{align*}
The proof of the other part is similar.
\end{proof}

Recall the braiding bicharacter $R=r\inv\ast r^{\tau}$ associated to
$r$, and write $S=S^{(R)}$.
\begin{prop}\label{prop:mrbraided}
  For any bicharacter $r$ on $M$ the twist $(M,m_r,1,S)$ is a braided
  algebra.
\end{prop}

\begin{proof} We need to check the axioms in Definition
  \ref{defn:Braidedalg}.  The vacuum and Yang-Baxter axioms are dealt
  with in Lemma \ref{lem:SrhoVacYangBaxter}. For unitarity we have
  $\tau\circ S\circ \tau=S^{(R^\tau)}$ so that by Lemma
  \ref{lem:Srhosigma} and \eqref{eq:AppRunitary}
\[
S\circ \tau \circ S\circ \tau=S^{(R)}\circ S^{(R^\tau)}=S^{(R\ast
  R^\tau)}=S^{(\epsilon)}=1_{M^{\otimes 2}}.
\]
Compatibility of $S$ with the multiplication $m_r$ is the case
$\rho=R$ of Lemma \ref{lem:rhocompatmr}.

Now $m_r$ is braided commutative:
\[
m_rS=m\circ S^{(r)}\circ S^{(R)}=m\circ S^{(r^\tau)}=m\circ \tau\circ
S^{(r)}\circ\tau=m_r\circ\tau,
\]
by Lemma \ref{lem:Srhosigma}, \eqref{eq:Appbraidingforr} and the fact
that $m$ is commutative. From compatibility of $S$ with multiplication
$m_r$ and the  Yang-Baxter equation it follows that $m_r$ is associative.  The
braiding axiom for $m_{r,3}=m_r(1\otimes m_r)=m_r(m_r\otimes1)$
follows from this.
\end{proof}

\subsection{Bicharacters and Group Action}
\label{subsec:BichGroup}
Now we assume that we have an action of a group $G$ on the commutative
and cocommutative Hopf algebra $M$ compatible with the multiplication
and the comultiplication: 
\[
gm=m\circ \Delta(g), \quad \Delta(gm)=\Delta_G(g)\Delta(m).
\]
Define for any bicharacter $r$ on $M$ and $g\in G$
\[
r^g=r\circ \Delta(g).
\]
It is easy to check that $r^g$ is again a bicharacter, so that we can
write
\begin{equation}
  \label{eq:defrgRg}
r^g=r\ast R^g,\quad R^g= r\inv\ast r^g.
\end{equation}
Also $R^g$ is then a bicharacter. Define $S^g=S^{(R^g)}$.
\begin{lem}\label{lem:Gsym} For all $g\in G$ and bicharacters $r$ on $M$
  \[
gm_r S^g=m_r \Delta(g).
\]
\end{lem}

\begin{proof} By Lemma \ref{lem:Srhosigma} and \eqref{eq:Appbraidingforr}
  \begin{align*}
    gm_rS^g&=gm\circ S^{(r)}\circ S^{(R^g)}=g m\circ S^{(r\ast
    R^g)}=\\
  &=gm\circ S^{(r^g)}=m\circ \Delta(g)\circ S^{(r^g)}=m\circ S^{(r)}\Delta(g).
  \end{align*}
Here we use
\begin{equation}
  \label{eq:DeltaSrg}
  \Delta(g)S^{(r^g)}=S^{(r)}\Delta(g),
\end{equation}
which follows from the definition of $S^{(r)}$, see \eqref{eq:DefSrho}.
\end{proof}

\begin{cor}
  Let $r$ be a bicharacter on a commutative and cocommutative Hopf
  algebra $M$ with an action of a group $G$. Then $(M,m_r,1,S)$ is a
  braided $G$-algebra for the maps
\[
S^g=S^{(R^g)},\quad g\in G,
\]
where $R^g$ is defined in \eqref{eq:defrgRg}.
\end{cor}
\begin{proof}We need to check the axioms in Definition
  \ref{defn:BraidedGalg}.  The $G$-symmetry axiom is verified in the
  previous Lemma \ref{lem:Gsym}. The Vacuum Axiom and $G$-Yang-Baxter
  Axiom for $S^g$ are verified in Lemma \ref{lem:SrhoVacYangBaxter},
  as $S^g=S^{(R^g)}$ and $R^g$ is a bicharacter. The compatibility of
  $S^g$ with multiplication $m_r$ is the case $\rho=R^g$ of Lemma
  \ref{lem:rhocompatmr}. For multiplicativity
  \begin{align*}
    S^{gh}&=S^{(r\inv\ast r^{gh})}=S^{(r\inv)}\circ S^{(r^{gh})}=\\
          &=S^{(r\inv)}\circ\Delta(gh)\inv\circ
          S^{(r)}\circ\Delta(gh)= &&\text{(by
          \eqref{eq:DeltaSrg})}\\
          &=S^{(r\inv)}\circ\Delta(h)\inv\circ S^{(r)}\circ S^{(r\inv)}\circ S^{(r^g)}\circ\Delta(h)=\\
          &=S^{(r\inv)}\circ  S^{(r^h)}\circ \Delta(h)\inv\circ
          S^{g}\circ\Delta(h)=\\
          &=S^h\circ\Delta(h\inv)\circ S^g\circ \Delta(h).
  \end{align*}
\end{proof}

\begin{remark}
  In a braided $G$-algebra we implement the action of $G$ by a system
  of maps $S^g$ satisfying
\[
gmS^g=m\Delta(g).
\]
In the bicharcter case of a twisted multiplication $m_r=m\circ S^{(r)}$ we can also implement
the group action by twisting the coproduct on $G$: we have
\begin{align*}
  g m_r&=g m\circ S^{(r)}=m\circ\Delta(g)\circ S^{(r)}=\\
       &=m\circ S^{(r)} \circ S^{(r\inv)}\circ\Delta(g)\circ
       S^{(r)}=\\
       &=m_r \Delta_r(g),
\end{align*}
where the twisted coproduct is
\[\Delta_r(g)=S^{(r\inv)}\circ\Delta(g)\circ
       S^{(r)}.
\]
The fact that the two approaches are equivalent, 
\[
\Delta(g)
(S^g)\inv=\Delta_r(g),
\]
follows from \eqref{eq:DeltaSrg}. It is at this point not clear
whether one can replace in an arbitrary braided $G$-algebra the maps $S^g$ by a
twist of the coproduct.
\end{remark}

\section{Braiding maps}
\label{app:braiding}

Let $V$ be a free $k$-module and let $\MapzV$ be the space of linear maps
\[
V^{\ten n}\to V^{\ten n}[z_i^{\pm 1},(z_i-z_j)\inv][[t]],\quad 1\le
i<j\le n.
\]
Suppose we are given $S_{z_1,z_2}\in\Map_{z_1,z_2}(V^{\ten 2})$ that
satisfies
\begin{align}\label{App;eq:propertiesS1}
S_{z_1,z_2}\circ\tau\circ S_{z_2,z_1}\circ \tau&=1_{V^{\ten 2}},\\
S_{z_1,z_2}^{12}S_{z_1,z_3}^{13}S_{z_2,z_3}^{23}&=S_{z_1,z_3}^{14}S_{z_1,z_3}^{13}S_{z_1,z_2}^{12}.
\label{App;eq:propertiesS2}
\end{align}
We then define for each $\mathbf{f}\in\mathcal{S}_n$ an element
$S^{\mathbf{f}}\subz\in\MapzV$ as follows. First, for
$\w_i\in\mathcal{S}_n$ a simple transposition, define
\[
S^{\w_i}\subz=1^{i-1}\ten S^{(\tau)}_{z_i,z_{i+1}}\ten i^{n-i-1},
\]
and extend this to $\mathbf{f}\in\mathcal{S}_n$ by expanding
it in simple transpositions and using
\begin{equation}\label{app;eq:defSw}
S^{\mathbf{f}{\mathbf{g}}}\subz=S^{\mathbf{g}}\subz \sigma_{\mathbf{g}}
S^\mathbf{f}_{{\mathbf{g}}\inv(z_1,\dots,z_n)}(\sigma_\mathbf{f})^{-1}.
\end{equation}
The problem is that the expansion of $\mathbf{f}$ is not unique,
because of the relations \eqref{eq:relSn} and \eqref{eq:relSntwo} in
$\mathcal{S}_n$.

To address this problem introduce the free monoid $\mathcal{F}_n$
generated by symbols $\tilde \w_i$, $i=1,2,\dots,n-1$. In
$\mathcal{F}_n$ any element $\tilde {\mathbf{f}}$ has a unique
expression in terms of the $\tilde \w_i$s. Consider the semi-direct
product $\MapzV\rtimes \mathcal{S}_n$: elements of the semi-direct
product are pairs $(A\subz,\mathbf{f})$, with product
\begin{equation}
  \label{App;eq:semidirectprod}
  (A\subz,\mathbf{f}).(B\subz,{\mathbf{g}})=(A\subz\circ \sigma_\mathbf{f}\inv\circ
B_{\mathbf{f}(z_1,z_2,\dots,z_n)}\circ \sigma_\mathbf{f}, \mathbf{f}{\mathbf{g}}).
\end{equation}
We have a homomorphism $\mathcal{F}_n\to \mathcal{S}_n$, which maps
generator $\tilde \w_i$ to simple transposition $\w_i$.  Let
\[
\phi\colon\mathcal{F}_n\to \MapzV\rtimes \mathcal{S}_n
\]
be given on generators by
\[
\phi(\tilde \w_i)=(S\subz^{\w_i},\w_i),
\]
and we extend this to all of $\mathcal{F}_n$ as an
\emph{anti-homomorphism} of monoids.

We need some more notation. If $\tilde {\mathbf{f}}=\tilde \w_{i_1}\tilde
\w_{i_2}\dots \tilde \w_{i_k}\in \mathcal{F}_n$, and the corresponding
permutation is $ \mathbf{f}= \w_{i_1}\w_{i_2}\dots \w_{i_k}\in\mathcal{S}_n$, then
introduce
\[
{\mathbf{g}}_\ell=\w_{i_k}\w_{i_{k-1}}\dots \w_{i_{\ell+1}},\quad \ell= 1, 2, \dots,
k-1,
\]
and ${\mathbf{g}}_k=1$.

\begin{lem}\label{lem:phiproperties}
  Let $\tilde {\mathbf{f}}=\tilde \w_{i_1}\tilde \w_{i_2}\dots \tilde \w_{i_k}$ and $\mathbf{f}=
\w_{i_1}\w_{i_2}\dots  \w_{i_k}$. Then
\[
\phi(\tilde {\mathbf{f}})=(S^\mathbf{f}\subz,\mathbf{f}\inv)\in
\MapzV\rtimes \mathcal{S}_n,
\]
where
\[
S^\mathbf{f}\subz=S^k S^{k-1}\dots S^1 \sigma_\mathbf{f}\inv, \quad
S^\ell=S^{\w_{i_\ell}}_{{\mathbf{g}}_\ell(z_1,z_2,\dots,z_n)}\tau_{i_\ell}.
\]
Furthermore, for $\tilde {\mathbf{f}},\tilde {\mathbf{g}}\in \mathcal{F}_n$
\begin{equation}
  \label{eq:multphifg}
  S^{\mathbf{f}{\mathbf{g}}}\subz\sigma_{\mathbf{f}{\mathbf{g}}}
  =S^{\mathbf{g}}\subz \sigma_{\mathbf{g}} S^\mathbf{f}_{{\mathbf{g}}\inv(z_1,z_2,\dots,z_n)}\sigma_\mathbf{f}.
\end{equation}
\end{lem}

\begin{proof}
  Using the anti-homomorphism property of $\phi$ and the
  multiplication \eqref{App;eq:semidirectprod} we have
\begin{align*} 
  \phi(\tilde {\mathbf{f}}) &=\phi(\tilde \w_{i_k})\phi(\tilde
  \w_{i_{k-1}})\dots
  \phi(\tilde \w_{i_1})=\\
  &=(S^{\w_{i_k}}\subz,\w_{i_k}).(S^{\w_{i_{k-1}}}\subz,\w_{i_{k-1}}).
  \dots.(S^{\w_{i_1}}\subz,\w_{i_1})=\\
  &=(S^{\w_{i_k}}\subz\tau_{i_k}S^{\w_{i_{k-1}}}_{\w_{i_k}(z_1,\dots,z_n)}\tau_{i_{k}},
  \w_{i_k}\w_{i_{k-1}}).(S^{\w_{i_{k-2}}}\subz,\w_{i_{k-2}}).\dots\\
  &\quad \dots(S^{\w_{i_1}}\subz,\w_{i_1})=\\
  &=(S^k S^{k-1}\dots S^1\sigma_\mathbf{f}\inv, \mathbf{f}\inv).
\end{align*}
  This proves the first part. Then
\begin{align*}
  \phi(\tilde {\mathbf{f}}\tilde {\mathbf{g}}) &=\phi(\tilde {\mathbf{g}})\phi(\tilde {\mathbf{f}})=\\
  &=(S^{{\mathbf{g}}}\subz,{\mathbf{g}}\inv).(S^{\mathbf{f}}\subz,\mathbf{f}\inv)=(S^{\mathbf{g}}\subz
  \sigma_{\mathbf{g}}
  S^\mathbf{f}_{{\mathbf{g}}\inv(z_1,\dots,z_n)}\sigma_{\mathbf{g}}\inv,
  {\mathbf{g}}\inv \mathbf{f}\inv).
\end{align*}
Since $\sigma_{\mathbf{g}}\inv \sigma_{\mathbf{f}{\mathbf{g}}}=\sigma_\mathbf{f}$ \eqref{eq:multphifg} follows.
\end{proof}
The observant reader might object to the notation
$S^{\mathbf{f}}\subz$ used in the above Lemma: this map depends a
priori on the element $\tilde{\mathbf{f}}\in \mathcal{F}_n$, not just
on its image in $\mathcal{S}_n$. The following Lemma justifies the
notation.

\begin{lem}\label{lem:phifactors}
  The map $\phi\colon \mathcal{F}_n\to \MapzV\rtimes \mathcal{S}_n$
  factors through the canonical map $\mathcal{F}_n\to\mathcal{S}_n$.
\end{lem}

\begin{proof}
  We need to check that the relations \eqref{eq:relSn} and
  \eqref{eq:relSntwo} (with $\w_i$ replaced by $\tilde \w_i$) belong to
  the kernel of $\phi$. But we have
  \begin{align*}
    \phi(\tilde \w_i^2)&=\phi(\tilde \w_i)\phi(\tilde \w_i)=\\
    &=(S^{\w_i}\subz,\w_i).(S^{\w_i}\subz,\w_i)=(S^{\w_i}\subz \tau_i
    S^{\w_i}_{\w_i(z_1,\dots,z_n)}\tau_i,\w_i\w_i)=\\
    &=(1,1)
  \end{align*}
by the definition \eqref{app;eq:defSw} and the property
\eqref{eq:propertiesS1}. Next, by the definition \eqref{app;eq:defSw} and
\eqref{eq:relSntwo} we have, if $\abs{i-j}\ge2$,
\[
\phi(\tilde \w_i \tilde \w_j)=\phi(\tilde \w_j)\phi(\tilde
\w_i)=\phi(\tilde \w_i)\phi(\tilde \w_j)=\phi(\tilde \w_j\tilde \w_i)
\]
Finally 
\[
\phi(\tilde \w_i\tilde \w_{i+1}\tilde \w_i)=\phi(\tilde \w_{i+1}\tilde
\w_{i}\tilde \w_{i+1})
\]
follows from the Yang-Baxter equation \eqref{eq:propertiesS2}.
\end{proof}

The conclusion is that Definition \ref{defn:braiding-maps-n2} of
$S^{\mathbf{f}}\subz$ is well defined.

\bibliographystyle{amsalpha}

\def\cprime{$'$}
\providecommand{\bysame}{\leavevmode\hbox to3em{\hrulefill}\thinspace}
\providecommand{\MR}{\relax\ifhmode\unskip\space\fi MR }
\providecommand{\MRhref}[2]{%
  \href{http://www.ams.org/mathscinet-getitem?mr=#1}{#2}
}
\providecommand{\href}[2]{#2}

\end{document}